\documentclass{article}
\usepackage{fullpage}

\newtheorem{theorem}{Theorem}
\newtheorem{lemma}{Lemma}
\newtheorem{proposition}{Proposition}

\begin{document}
\title{On Classification of Normal Operators \\
in Real Spaces with Indefinite Scalar Product}
\author{O.V.Holtz \kern0.07\textwidth V.A.Strauss \\
\small Department of Applied Mathematics \\
\small Chelyabinsk State Technical University \\
\small 454080 Chelyabinsk, Russia}
\date{}
\maketitle

\begin{abstract}
A real finite dimensional space with indefinite scalar product having
$v_{-}$ negative squares and $v_{+}$ positive ones is considered.  The
paper presents a classification of operators that are normal with
respect to this product for the cases $min\{v_{-},v_{+}\}=1, \;
2$. The approach to be used here was developed in the papers \cite{1}
and \cite{2}, where the similar classification was obtained for
complex spaces with $v=min\{v_{-},v_{+}\}=1, \; 2$, respectively.
\end{abstract}

\section{Introduction}
Consider a real linear space $R^{n}$ with an indefinite scalar product
$[ \cdot \: , \cdot ]$.  By definition, the latter is a nondegenerate
sesquilinear Hermitian form.  If the ordinary scalar product $(\cdot
\:, \cdot)$ is fixed, then there exists a nondegenerate Hermitian
operator $H$ such that $[x,y]=(Hx,y)$ $\forall x,y\in R^{n}$.  If $A$
is a linear operator ($A: R^{n} \rightarrow R^{n}$), then the
{\em{$H$-adjoint\/}} of $A$ (denoted by $A^{[*]}$) is defined by the
identity $[A^{[*]}x,y] \equiv [x,Ay]$. An operator $N$ is called
{\em{$H$-normal\/}} if $NN^{[*]}=N^{[*]}N$, an operator $U$ is called
{\em{$H$-unitary\/}} if $UU^{[*]}=I$, where $I$ is the identity
transformation.

Let $V$ be a nontrivial subspace of $R^{n}$. The subspace $V$ is
called {\em{neutral}} if $[x,y]=0$ $\forall x,y \in V$. If from the
conditions $x \in V$ and $\forall y \in V \; [x,y]=0$ it follows that
$x=0$, then $V$ is called {\em{nondegenerate}}.  The subspace
$V^{[\perp]}$ is defined as the set of all vectors $x \in R^{n}$:
$[x,y]=0$ $\forall y \in V$. If $V$ is nondegenerate, then
$V^{[\perp]}$ is also nondegenerate and $V \dot{+} V^{[\perp]}=R^{n}$.

A linear operator $A$ acting in $R^{n}$ is called {\em{decomposable}}
if there exists a nondegenerate subspace $V \subset R^{n}$ such that
both $V$ and $V^{[\perp]}$ are invariant for $A$ or (it is the same)
if $V$ is invariant both for $A$ and $A^{[*]}$. Then $A$ is the
{\em{orthogonal sum of $A_{1}=A|_{V}$ and
$A_{2}=A|_{V^{[\perp]}}$}}. If an operator $A$ is not decomposable, it
is called {\em{indecomposable}}.

Throughout what follows by a rank of a space we mean
$v=\min\{v_{-},$$v_{+}\}$, where $v_{-}\:$($v_{+}$) is the number of
negative (positive) squares of the quadratic form $[x,x]$, i.e., the
number of negative (positive) eigenvalues of the operator $H$.  Note
that without loss of generality it can be assumed that $v_{-} \leq
v_{+}$ (otherwise $H$ can be replaced by $-H$; the latter
(nondegenerate Hermitian operator) has opposite eigenvalues). Later on
we assume that $v_{-} \leq v_{+}$.

The problem is to obtain a complete classification for $H$-normal
operators acting in $R^{n}$, i.e., to find a set of canonical forms
such that any $H$-normal operator could be reduced to one and only one
of these forms. Since it is sufficient to solve the problem only for
indecomposable operators, for any nondegenerate Hermitian matrix $H$
and for any indecomposable $H$-normal matrix $N$ we would like to
point out one and only one of the canonical pairs of matrices $\{
\tilde{N}, \tilde{H} \}$ so that the pair $\{N,H\}$ is unitarily
similar to $\{ \tilde{N}, \tilde{H} \}$ (two pairs of matrices
$\{N_{1},H_{1}\}$ and $\{N_{2},H_{2}\}$, where $H_{1}$ and $H_{2}$ are
nondegenerate Hermitian matrices, are called {\em{unitarily similar}}
if $N_{2}=T^{-1}N_{1}T$, $H_{2}=T^{*}H_{1}T$ for some invertible
matrix $T$; if $H_{1}=H_{2}$, then they are {\em{$H_{1}$-unitarily
similar}}). In what follows such a classification is presented for
operators acting in spaces of rank $1$ and $2$.  As in \cite{2}, we
will denote by $I_{r}$ the identity matrix of order $r \times r$, by
$D_{r}$ the $r \times r$ matrix with 1's on the secondary diagonal and
zeros elsewhere, and by $A \oplus B \oplus \ldots \oplus C$ a block
diagonal matrix with blocks $A$, $B$, $\ldots$, $C$.

We are grateful to Prof. Leiba Rodman for his attention to our work
and very helpful comments on this paper.

\section{On Decomposition of $H$-normal Operators in Real Spaces}

Let an $H$-normal operator $N$ act in $R^{n}$ and have $p$ distinct
real eigenvalues $\lambda_{1}, \lambda_{2}, \ldots, \lambda_{p}$ and
$q$ distinct pairs of complex conjugate eigenvalues $\alpha_{p+1}\pm
i\beta_{p+1}, \alpha_{p+2}\pm i\beta_{p+2}, \ldots, \alpha_{p+q}\pm
i\beta_{p+q}$.  Let us define $$ \varphi_{k}(\lambda)=\left\{
\begin{array}{ccc} (\lambda-\lambda_{k})^{n}, & \mbox{if} & 1 \leq k
\leq p \\
(\lambda^{2}-2\alpha_{k}\lambda+\alpha_{k}^{2}+\beta_{k}^{2})^{n}, &
\mbox{if} & p < k \leq p+q, \end{array} \right. $$ $$Q_{ij}=\{x: \;
\varphi_{i}(N)x=\varphi_{j}(N^{[*]})x=0\}, \; \;i,j=1,\ldots,p+q,$$
$$\Omega=\{(i,j): \; Q_{ij}\neq \{0\} \}.$$

\begin{proposition}

  The subspaces $Q_{ij}$ have the following properties:
(1) $Q_{ij} \cap Q_{kl}=\{0\}$ $\; \forall \; (i,j) \neq (k,l)$.

(2) $\sum_{(i,j) \in \Omega} Q_{ij}=R^{n}$.

(3) Each subspace $Q_{ij}$ is invariant for both $N$ and $N^{[*]}$.
(4) Eigenvalues of the operator $N|_{Q_{ij}}$ are roots of
    $\varphi_{i}(\lambda)$, those of the operator $N^{[*]}|_{Q_{ij}}$
    are roots of $\varphi_{j}(\lambda)$.
(5) $[Q_{ij},Q_{kl}]=0$  $\; \forall \; (i,j) \neq (l,k)$.
\end{proposition}

{\bf{Proof:}}

 (1) Suppose $(i,j) \neq (k,l)$. Without loss of generality it can be
     assumed that $i \neq k$. Let $\exists x:$ $x \in Q_{ij}$, $x \in
     Q_{kl}$, i.e., $\varphi_{i}(N)x=\varphi_{k}(N)x=0$. Since the
     polynomials $\varphi_{i}(\lambda)$ and $\varphi_{k}(\lambda)$ are
     relatively prime, there exist polynomials $\psi_{i}(\lambda)$,
     $\psi_{k}(\lambda)$ such that the matrix identity $I \equiv
     \psi_{i}(A)\varphi_{i}(A)+\psi_{k}(A)\varphi_{k}(A)$ is
     valid. Consequently,
     $x=\psi_{i}(N)\varphi_{i}(N)x+\psi_{k}(N)\varphi_{k}(N)x=0$.

(2) The greatest common divisor of the polynomials
    $\xi_{1}(\lambda)=\prod_{i \neq 1}\varphi_{i}(\lambda)$,
    $\xi_{2}(\lambda)=\prod_{i \neq 2}\varphi_{i}(\lambda)$, $\ldots$,
    $\xi_{p+q}(\lambda)=\prod_{i \neq p+q}\varphi_{i}(\lambda)$ is
    equal to $1$, therefore, there exist polinomials
    $\psi_{1}(\lambda)$, $\psi_{2}(\lambda)$, $\ldots$,
    $\psi_{p+q}(\lambda)$ such that $I=\sum_{i=1}^{p+q}
    \psi_{i}(A)\xi_{i}(A)$ $\;\forall A$.  Hence, $\forall x$
    $x=\sum_{i=1}^{p+q} \psi_{i}(N)\xi_{i}(N)x=\sum_{i=1}^{p+q}x_{i}$
    (where $x_{i}=\psi_{i}(N)\xi_{i}(N)x$).  Since the product of all
    $\varphi_{i}(\lambda)$ annihilates $N$, we have
    $\varphi_{i}(N)x_{i}=0$ $\forall i$, i.e.,
    $R^{n}=\sum_{i=1}^{p+q}Q_{i}$, where $Q_{i}=\{x: \;
    \varphi_{i}(x)=0\}$.  Similarly, each subspace $Q_{i}$ is a direct
    sum of the subspaces $Q_{ij}=\{x \in Q_{i}: \;
    \varphi_{j}(N^{[*]})x=0\}$.  Disregarding the trivial subspaces
    $Q_{ij}$, we obtain the desired equality $R^{n}=\sum_{(i,j)\in
    \Omega} Q_{ij}$.

(3) Since $N$ and $N^{[*]}$ commute, for all $(i,j)$ and $x \in
    Q_{ij}$ we have $0=N\varphi_{i}(N)x=\varphi_{i}(N)Nx$,
    $0=N\varphi_{j}(N^{[*]})x=\varphi_{j}(N^{[*]})Nx$, i.e., $Nx \in
    Q_{ij}$. It can be checked in the same way that $N^{[*]}x \in
    Q_{ij}$.

(4) Let $N|_{Q_{ij}}$ have an eigenvalue $\lambda_{0}$ such that
    $\varphi_{i}(\lambda_{0}) \neq 0$. Then there exists a (real or
    complex) eigenvector $x \neq 0$ corresponding to the eigenvalue
    $\lambda_{0}$.  Since the polynomials $\lambda-\lambda_{0}$ and
    $\varphi_{i}(\lambda)$ are relatively prime, there exist
    polynomials $\psi_{1}(\lambda)$, $\psi_{2}(\lambda)$ such that the
    identity $I=\psi_{1}(A)(A-\lambda_{0}I)+\psi_{2}(A)\varphi_{i}(A)$
    holds for all (complex) matrices $A$.  Consequently,
    $x=\psi_{1}(N)(N-\lambda_{0}I)x+\psi_{2}(N)\varphi_{i}(N)x=0$
    because $(N-\lambda_{0}I)x=\varphi_{i}(N)x=0$. The contradiction
    obtained shows that all eigenvalues of $N|_{Q_{ij}}$ are roots of
    $\varphi_{i}(\lambda)$. The operator $N^{[*]}|_{Q_{ij}}$ can be
    considered in the same way.

(5) Let $i \neq l$. Take arbitrary vectors $x \in Q_{ij}$, $y \in
    Q_{kl}$.  Since the eigenvalues of $N|_{Q_{ij}}$ are not roots of
    $\varphi_{l}(\lambda)$, the operator $\varphi_{l}(N)|_{Q_{ij}}$ is
    nondegenerate. Therefore, $\exists z \in Q_{ij}:$
    $\varphi_{l}(N)z=x$.  We have
    $[x,y]=[\varphi_{l}(N)z,y]=[z,\varphi_{l}(N^{[*]})y]=[z,0]=0$.

The proof of the proposition is completed.

Now let $V_{i}=Q_{ii}$ $((i,i) \in \Omega)$,
$V_{jk}=span\{Q_{jk},Q_{kj}\}$ ($(j,k) \in \Omega, \; j<k$). The
subspaces $V_{i}$, $V_{jk}$ are mutually orthogonal, the intersection
of any two of them is zero, and their sum is $R^{n}$. It follows from
the nondegeneracy of $H$ that each subspace $V_{i}$, $V_{jk}$ is
nondegenerate.  The restriction $N|_{V_{i}}$ has the only real
eigenvalue $\lambda_{i}$ if $i \leq p$ or the pair of complex
conjugate eigenvalues $\alpha_{i} \pm i\beta_{i}$ if $i>p$. The
restriction $N|_{V_{jk}}$ has two distinct real eigenvalues
$\lambda_{j}$, $\lambda_{k}$ if $j,k \leq p$, one real eigenvalue
$\lambda_{j}$ and the pair of complex conjugate eigenvalues
$\alpha_{k}\pm i\beta_{k}$ if $j \leq p$, $k>p$, or two distinct pairs
$\alpha_{j} \pm i\beta_{j}$, $\alpha_{k} \pm i\beta_{k}$ if $j,k>p$.

Thus, we have proved the following lemma:

\begin{lemma}

  Any $H$-normal operator $N$ acting in $R^{n}$ is an orthogonal sum
  of $H$-normal operators each of which has one of the following sets
  of eigenvalues: 
\begin{description}
\item{(a)} one real eigenvalue; \item{(b)} two distinct real
eigenvalues; \item{(c)} two complex conjugate eigenvalues; \item{(d)}
one real and two complex conjugate eigenvalues; \item{(e)} two
distinct pairs of complex conjugate eigenvalues.  \end{description}

\end{lemma}

This lemma shows the principal difference between real and complex
spaces because indecomposable operators acting in complex spaces have
either one or two distinct eigenvalues (Lemma~1 from~\cite{1}).

\section{Classification of $H$-normal Operators Acting in Spaces of Rank $1$}

This section is closely related to \cite{1}.

Let us classify indecomposable $H$-normal operators acting in a space
$R^{n}$ of rank $1$. According to Lemma~1, we can consider only
operators having one of the sets of eigenvalues (a) - (e). However,
for a space of rank $1$ not all variants are possible, namely, the
alternatives (d) and (e) cannot be realized.  Indeed, if $N|_{Q_{12}}$
(or $N^{[*]}|_{Q_{12}}$) has two eigenvalues $\alpha \pm i \beta$, the
subspace $Q_{12}$ is necessarily of dimension $2$ or higher. However,
since $Q_{12}$ is neutral, $dim \: Q_{12} \leq 1$.  Thus, the
alternatives (d) and (e) are impossible.  Let us consider the
remaining variants and prove the following theorem:

\begin{theorem}

If an indecomposable $H$-normal operator $N$ ($N: \: R^{n} \rightarrow
R^{n}$) acts in a space with indefinite scalar product having
$v_{-}=1$ negative squares and $v_{+} \geq 1$ positive ones, then $2
\leq n \leq 4$ and the pair $\{N,H\}$ is unitarily similar to one and
only one of canonical pairs (\ref{th1.1}), (\ref{th1.2}),
(\ref{th1.3}), (\ref{th1.4}), (\ref{th1.5}), (\ref{th1.6}):
\begin{equation} N=\left( \begin{array}{cc} \lambda_{1} & 0 \\ 0 &
\lambda_{2} \end{array} \right), \;\; \lambda_{1}< \lambda_{2}, \;\;
H=D_{2}, \label{th1.1} \end{equation}

\begin{equation}
N=\left( \begin{array}{cc} \alpha & \beta \\ -\beta & \alpha
\end{array} \right), \;\; \beta>0, \;\; H=D_{2}, \label{th1.2}
\end{equation} \begin{equation} N=\left( \begin{array}{cc} \lambda & z
\\ 0 & \lambda \end{array} \right), \;\; z=\pm 1, \;\; H=D_{2},
\label{th1.3} \end{equation}
\begin{equation} N=\left( \begin{array}{ccc} \lambda & 1 & 0   \\ 
0 & \lambda & 1  \\ 0 & 0 & \lambda \end{array} \right), \;\; H=D_{3}, \label{th1.4} \end{equation} 
\begin{equation} N=\left( \begin{array}{ccc} \lambda & 1 & r   \\ 
0 & \lambda & -1  \\ 0 & 0 & \lambda \end{array} \right), \;\; H=D_{3}, \label{th1.5} \end{equation} 
\begin{equation} N=\left( \begin{array}{cccc} \lambda  & 1 & 0 & 0  \\ 0 & \lambda & 0 & \cos \alpha 
 \\ 0 & 0 & \lambda & \sin  \alpha  \\ 0 & 0 & 0 & \lambda \end{array} \right), \; 0<\alpha<\pi, \;\; 
H=\left( \begin{array}{ccc} 0 & 0 & 1  \\ 0 & I_{2} & 0  \\ 1 & 0 & 0  \end{array} \right). \label{th1.6} \end{equation} \end{theorem} 

The proof of the theorem is presented in the following subsections.
\subsection{One Real Eigenvalue of $N$} Let us take advantage of
Proposition~1 from~\cite{2}, which is proved for complex spaces but is
valid for real ones as well: {\em{If an indecomposable $H$-normal
operator $N: \;R^{n} \rightarrow R^{n}$ ($n>1$) has the only eigenvalue
$\lambda$, then there exists a decomposition of $R^{n}$ into a direct sum of subspaces \begin{equation} S_{0}=\{x: \; (N-\lambda I)x=(N^{[*]}- \lambda I)x=0\}, \label{aug1} \end{equation} $S$, $S_{1}$ such that

\begin{equation} N=\left( \begin{array}{ccc} N'=\lambda I & * & * \\ 0 & N_{1} & * \\ 0 & 0 & N''=\lambda I \end{array} \right), \;\; H=\left( \begin{array}{ccc} 0 & 0 & I \\ 0 & H_{1} & 0 \\ I & 0 & 0 \end{array} \right), \label{pred1} \end{equation} where $N':$ $S_{0} \rightarrow S_{0}$, $N_{1}:$ $S \rightarrow S$, $N'':$ $S_{1} \rightarrow S_{1}$, the internal operator $N_{1}$ is $H_{1}$-normal, and the pair $\{N_{1},H_{1}\}$ is determined up to unitary similarity. To go over from one decomposition $R^{n}=S_{0}\dot{+}S\dot{+}S_{1}$ to  another by a transformation $T$ it is necessary that the matrix $T$ be block triangular with respect to both decompositions.}}

Since $S_{0}$ is neutral, $dim \: S_{0}=1$. According to Proposition~2
from~\cite{2}, if the subspace $S_{0}$ is one-dimensional, then the
operator $N$ is indecomposable. So, it is not necessary to check the
indecomposability for each canonical form to be obtained in this
subsection.  As $H$ has one negative eigenvalue, $H_{1}$ has only
positive eigenvalues and one can assume that $H_{1}=I$, $N_{1}=
\lambda I$. Later on we will no longer stipulate that $H_{1}=I$,
$N_{1}=\lambda I$. By Theorem~1 of ~\cite{2} (it is also valid for
real spaces), $n \leq 4$. Consider the cases $n=2$, $3$, $4$
successively.

\subsubsection{$n=2$}
The matrices $N$ and $H$ have form (\ref{pred1}): $$N=\left
( \begin{array}{cc} \lambda & a \\ 0 & \lambda \end{array} \right),
\;\; H=\left( \begin{array}{cc} 0 & 1 \\ 1 & 0 \end{array} \right).$$
Since $S_{0} \cap S_{1} = \{0\}$, $a\neq 0$. Let
$\widetilde{v_{1}}=\sqrt{|a|}v_{1}$, $\widetilde{v_{2}}=1/\sqrt{|a|}
v_{2}$.  Then we do not change the matrix $H$ and reduce $N$ to form
(\ref{th1.3}).  Since (\ref{th1.3}) is a special case of canonical
form (16) from Theorem~1 (\cite{1}), the number $z$ is an $H$-unitary
invariant, i.e., two forms (\ref{th1.3}) with different values of $z$
are not $H$-unitarily similar.

\subsubsection{$n=3$}
The matrices $N$ and $H$ have form (\ref{pred1}): $$ N=\left
( \begin{array}{ccc} \lambda & a & b \\ 0 & \lambda & c \\ 0 & 0 &
\lambda \end{array} \right), \;\; H=\left( \begin{array}{ccc} 0 & 0 &
1 \\ 0 & 1 & 0 \\ 1 & 0 & 0 \end{array} \right). $$ The condition of
the $H$-normality of $N$ is
$$ a^{2}=c^{2}. $$ If $a=0$, then $c=0$ and $v_{2} \in S_{0}$, which is 
impossible because of the condition $S_{0} \cap S=\{0\}$.  Therefore, 
$a \neq 0$. Let $\widetilde{v_{1}}=av_{1}$, $\widetilde{v_{3}}=1/a \: v_{3}$. 
then we reduce $N$ to the form $$   N=\left( \begin{array}{ccc}
\lambda & 1 & b' \\ 0 & \lambda & x \\ 0 & 0 & \lambda \end{array}
\right), \;\; x=\pm 1 $$ without changing the matrix $H$. If $x=1$,
take the $H$-unitary transformation $T$ (throughout what follows only
$H$-unitary transformations are used unless otherwise stipulated):$$ 
T=\left( \begin{array}{ccc} 1 & \frac{1}{2}b' & -\frac{1}{8}b'^{2} \\
0 & 1 & -\frac{1}{2}b' \\ 0 & 0 & 1 \end{array} \right). $$ It reduces
$N$ to form (\ref{th1.4}).  If $x=-1$, the number $b'$ turns out to be
$H$-unitary invariant.  Indeed, let $$N- \lambda I=\left (
\begin{array}{ccc} 0 & 1 & r \\ 0 & 0 & -1 \\ 0 & 0 & 0 \end{array}
\right), \;\; \tilde{N} - \lambda I=\left( \begin{array}{ccc} 0 & 1 &
\tilde{r} \\ 0 & 0 & -1 \\ 0 & 0 & 0
\end{array} \right), $$ and some matrix $T=\{t_{ij}\}_{i,j=1}^{3}$ satisfy 
The conditions \begin{eqnarray} NT & = & T\tilde{N}, \label{similar} \\ 
tT^{[*]} & = & I. \label{unitar} \end{eqnarray}
Then, according to Proposition~1 from~\cite{2}, $T$ is block
triangular with respect to the decomposition $R^{n}=S_{0}\dot{+}S
\dot{+} S_{1}$, i.e., upper triangular. Condition (\ref{similar})
implies \begin{eqnarray} & t_{11}=t_{22}=t_{33}, & \nonumber \\ &
t_{23}+rt_{33}=\tilde{r}t_{11}-t_{12}. & \label{smain}
\end{eqnarray} Since the diagonal terms of $T$ are equal to each other, 
From (\ref{unitar}) it follows that $t_{12}+t_{23}=0$.  Then from 
(\ref{smain}) we get $r=\tilde{r}$, Q.E.D.  The forms obtained are not 
$H$-unitarily similar. Indeed, let an $H$-unitary matrix 
$T=\{t_{ij}\}_{i,j=1}^{3}$ reduce the first
form to the second. Since $T$ is upper triangular (Proposition~1
from~\cite{2}), from (\ref{similar}) it follows that
$t_{11}=t_{22}=-t_{33}$, which is impossible because condition
(\ref{unitar}) implies $t_{11}t_{33}=1$.  Thus, we have obtained two
canonical forms: (\ref{th1.4}) and (\ref{th1.5}).

\subsubsection{$n=4$} 
The matrices $N$ and $H$ have form (\ref{pred1}): $$ N=\left
		 ( \begin{array}{cccc} \lambda & a & b & c \\ 0 &
		 \lambda & 0 & d \\ 0 & 0 & \lambda & e \\ 0 & 0 & 0 &
		 \lambda
\end{array} \right), \;\; 
H=\left( \begin{array}{cccc} 0 & 0 & 0 & 1 \\ 0 & 1 & 0 & 0 \\ 0 & 0 &
1 & 0 \\ 1 & 0 & 0 & 0
\end{array} \right). $$ 
The condition of the $H$-normality of $N$ is \begin{equation}
a^{2}+b^{2}=d^{2}+e^{2}. \label{h_norm} \end{equation} Since
$a^{2}+b^{2} \neq 0$ (otherwise $v_{2},v_{3} \in S_{0}$, which is
impossible), without loss of generality it can be assumed that $a \neq
0$.  Taking $\widetilde{v_{1}}=av_{1}$, $\widetilde{v_{4}}=v_{4}/a$,
we reduce $N$ to the form 
$$ N=\left (
\begin{array}{cccc} \lambda & 1 & b' & c' \\ 0 & \lambda & 0 & d' \\ 0
& 0 & \lambda & e' \\ 0 & 0 & 0 & \lambda \end{array} \right). $$
Further, let us apply the transformation $$ T=\left (
\begin{array}{cccc} \sqrt{1+b'^{2}} & 0 & 0 & 0 \\ 0 &
1/\sqrt{1+b'^{2}} & -b'/\sqrt{1+b'^{2}} & 0 \\ 0 & b'/\sqrt{1+b'^{2}}
& 1/\sqrt{1+b'^{2}} & 0 \\ 0 & 0 & 0 & 1/\sqrt{1+b'^{2}} \end{array}
\right). $$ Then we get
$$ N=\left( \begin{array}{cccc} \lambda  & 1 & 0 & c''  \\ 0 & \lambda & 0 & d'' 
 \\ 0 & 0 & \lambda & e''  \\ 0 & 0 & 0 & \lambda \end{array} \right). $$

Note that $e'' \neq 0$ because otherwise $v_{3} \in S_{0}$, which is
impossible because $S_{0} \cap S=\{0\}$.  The number $e''$ can be
replaced by $-e''$ by means of the ($H$-unitary) transformation
$\widetilde{v_{3}}=-v_{3}$. So, we can assume $e''>0$. Moreover, it
can be assumed that $c''=0$. To this end it is sufficient to take the
transformation
$$ T=\left( \begin{array}{cccc}
1 & 0 & c''/e'' & -\frac{1}{2}c''^{2}/e''^{2} \\ 0 & 1 & 0 & 0 \\ 0 &
0 & 1 & -c''/e'' \\ 0 & 0 & 0 & 1 \end{array} \right).$$ Then $c''$
will vanish, $d''$ and $e''$ will not change.  Condition
(\ref{h_norm}) of the $H$-normality of $N$ implies $d''=\cos \alpha$,
$e''=\sin \alpha$ ($\alpha \: \in \:(0; \pi)$).  Show the $H$-unitary
invariance of the parameter $\alpha$. Let an $H$-unitary matrix
$T=\{t_{ij}\}_{i,j=1}^{4}$ reduce $N$ to the form $$ \tilde{N}=\left (
\begin{array}{cccc} \lambda & 1 & 0 & 0 \\ 0 & \lambda & 0 & \cos
\tilde{\alpha} \\ 0 & 0 & \lambda & \sin \tilde{\alpha} \\ 0 & 0 & 0 &
\lambda \end{array} \right), \; \tilde{\alpha} \in (0;\pi).$$ Then,
according to Proposition~1 from~\cite{2}, $T$ is block triangular with
respect to the decomposition $R^{n}=S_{0}\dot{+}S\dot{+}S_{1}$ and
from (\ref{similar}) it follows that $t_{23}=0$. Now condition
(\ref{unitar}) yields $t_{32}=0$ . Applying (\ref{similar}) again, we
have \begin{eqnarray*} t_{11} & = & t_{22}, \\ t_{44} \cos \alpha & = &
t_{22} \cos \tilde{\alpha}, \\ t_{44} \sin \alpha & = & t_{33} \sin
\tilde{\alpha}.  \end{eqnarray*} Condition (\ref{unitar}) yields
$t_{11}t_{44}=t_{22}^{2}=t_{33}^{2}=1$ so that
$t_{11}=t_{22}=t_{44}=\pm 1$. Hence, $\cos \alpha=\cos
\tilde{\alpha}$. Since $\sin \alpha,\; \sin \tilde{\alpha}>0$, we have
$t_{33}=t_{44}$ and $\sin \alpha =\sin \tilde{\alpha}$. Consequently,
$\tilde{\alpha}=\alpha$, Q.E.D.  Thus, we have obtained canonical form
(\ref{th1.6}).

\subsection{Two Distinct Real Eigenvalues of $N$} According to Proposition~1, 
in this case $$ N=\left( \begin{array}{cc}
\lambda_{1} & 0 \\
                 0  & \lambda_{2}
	    \end{array} \right), \; \;
   H=\left( \begin{array}{cc}
                 0 & a \\
                 a & 0
            \end{array} \right), \; a \neq 0. $$
It can be assumed that $a=1$ (to this end it is sufficient to take
$\widetilde{v_{1}}=v_{1}/a$, $\widetilde{v_{2}}=v_{2}$). Since the order
of eigenvalues is not fixed, we can assume that $\lambda_{1} < \lambda_{2}$.
Thus, we have obtained canonical pair (\ref{th1.1}).

\subsection{Two Complex Conjugate Eigenvalues of $N$}
Let $N$ have two distinct eigenvalues $\lambda=\alpha + i \beta$,
$\overline{\lambda}=\alpha - i \beta$. Since $N$ and $N^{[*]}$
commute, there exists a vector $z=x+iy$ ($x,y \in R^{n}$) such
that either $Nz=\lambda z$, $N^{[*]}z=\overline{\lambda}z$ or
$Nz=\lambda z$, $N^{[*]}z=\lambda z$. In the first case $[z,\overline{z}]=0$.
Indeed,
$\overline{\lambda}[z,\overline{z}]=[\lambda z,\overline{z}]=[Nz,\overline{z}]=[z,N^{[*]}\overline{z}]=[z,\lambda\overline{z}]=\lambda[z,\overline{z}]$.
Therefore, $(\lambda- \overline{\lambda})[z,\overline{z}]=0$, hence
$[z,\overline{z}]=0$. Let us write in detail the condition obtained:
$[x+iy,x-iy]=[x,x]-i[y,x]-i[x,y]-[y,y]=0$, i.e., $[x,y]=0$, $[x,x]=[y,y]$.
Since two-dimensional subspace $V=span\{x,y\}$ cannot be neutral,
we have $[x,x] \neq 0$. Thus, $V$ is a nondegenerate subspace which is invariant for
$N$ and $N^{[*]}$. For $N$ to be indecomposable it is necessary to have
$R^{n}=V$. But $[x,x]=[y,y]$, i.e., $H$ is either positive or negative definite, which contradicts the condition $min \{ v_{-}, v_{+} \}=1$.
Thus, only the case $Nz=\lambda z$, $N^{[*]}z=\lambda z$ is possible.  It can be shown as before that $[z,z]=0$, i.e., $[x,x]=-[y,y]$ so that the subspace $V=span\{x,y\}$ is either nondegenerate
or neutral. As above, we see that $V$ is necessarily nondegenerate and
$V=R^{n}$.  

Thus, for the basis $\{x,y\}$ we have $$ N=\left( \begin{array}{cc} \alpha & \beta  \\ -\beta  & \alpha \end{array} \right), \;\; H=\left( \begin{array}{cc} a & b  \\
         b  & -a \end{array} \right)\; (a^{2}+b^{2} \neq 0). $$

Let us reduce $H$ to the form $D_{2}$ without changing the matrix $N$. To this end it is sufficient to take $$T=\left( \begin{array}{cc} t_{11} & t_{12} \\ -t_{12} & t_{11} \end{array} \right),$$
where \begin{eqnarray*} -2t_{11}t_{12} & = & a, \\ t_{11}^{2}-t_{12}^{2} & = & b \end{eqnarray*} (it can be checked that this system always has a real solution
$\{ t_{11},t_{12} \}$). Then
$$\left( \begin{array}{cc} a & b \\ b & -a \end{array} \right)= T^{*}\left( \begin{array}{cc} 0 & 1 \\
1 & 0 \end{array} \right)T, \;\;\; TN=NT. $$ One can replace $\beta$ by $-\beta$ by means of the $H$-unitary transformation $T=D_{2}$, therefore, one can assume that $\beta > 0$. Thus, we have obtained canonical pair (\ref{th1.2}). The proof of Theorem~1 is completed.

\section{Classification of $H$-normal Operators Acting in Spaces of Rank $2$}

The objective of this section is to prove the following theorem (the subspace $S_{0}$ and the internal operator $N_{1}$ are defined in Section 3.1 by formulas (\ref{aug1}), (\ref{pred1}), respectively):

\begin{theorem}
If an indecomposable $H$-normal operator $N$ ($N: \: R^{n} \rightarrow
R^{n}$) acts in a space with indefinite scalar product having
$v_{-}=2$ negative squares and $v_{+} \geq 2$ positive ones, then $4
\leq n \leq 8$ and the pair $\{N,H\}$ is unitarily similar to one and
only one of the canonical pairs \{(\ref{lemma1.1}),(\ref{lemma1.2})\}
- \{(\ref{lemma16.1}),(\ref{lemma16.2})\}. The list of all the
canonical pairs is as follows.

If $N$ has one real eigenvalue $\lambda$, $dim \:S_{0}=1$, the
internal operator $N_{1}$ is indecomposable, and $n=4$, then the pair
$\{N,H\}$ is unitarily similar to the canonical pair
\{(\ref{lemma1.1}),(\ref{lemma1.2})\}: \begin{equation} N=\left(
\begin{array}{cccc} \lambda & 1 & 0 & 0 \\ 0 & \lambda & z & 0 \\ 0 &
0 & \lambda & 1 \\ 0 & 0 & 0 & \lambda
\end{array} \right), \;\; z=\pm 1, \label{lemma1.1} \end{equation} 
\begin{equation}
H=D_{4}.  \label{lemma1.2} \end{equation} If $N$ has one real
eigenvalue $\lambda$, $dim \:S_{0}=1$, $N_{1}$ is indecomposable, and
$n=5$, then the pair $\{ N,H \}$ is unitarily similar to one and only
one of the canonical pairs \{(\ref{lemma2.1}),(\ref{lemma2.3})\},
\{(\ref{lemma2.2}),(\ref{lemma2.3})\}: \begin{equation} N=\left(
\begin{array}{ccccc} \lambda & 1 & 0 & 0 & 0 \\ 0 & \lambda & 1 & 0 &
0 \\ 0 & 0 & \lambda & 1 & 0 \\ 0 & 0 & 0 & \lambda & 1 \\ 0 & 0 & 0 &
0 & \lambda \end{array} \right), \label{lemma2.1} \end{equation} 
\begin{equation} N=\left( \begin{array}{ccccc} \lambda & 1 & -r_{1} & 0 
& r_{2} \\0 & \lambda & 1 & r_{1} & 0 \\ 0 & 0 & \lambda & -1 & -r_{1}
\\ 0 & 0 & 0 & \lambda & -1 \\ 0 & 0 & 0 & 0 & \lambda
\end{array} \right), \label{lemma2.2} \end{equation} \begin{equation}
H=D_{5}.  \label{lemma2.3} \end{equation} If $N$ has one real
eigenvalue $\lambda$, $dim\:S_{0}=1$, $N_{1}$ is decomposable, and
$n=4$, then the pair $\{ N,H \}$ is unitarily similar to one and only
one of the canonical pairs \{(\ref{lemma3.1}),(\ref{lemma3.3})\},
\{(\ref{lemma3.2}),(\ref{lemma3.3})\}: \begin{equation} N=\left
( \begin{array}{cccc} \lambda & 1 & 0 & 0 \\ 0 & \lambda & 0 & z \\ 0
& 0 & \lambda & 0 \\ 0 & 0 & 0 & \lambda \end{array} \right), \; z=
\pm 1, \label{lemma3.1} \end{equation}
\begin{equation} N=\left( \begin{array}{cccc} \lambda & 1 & z & 0 \\
 0 & \lambda & 0 & r \\ 0 & 0 & \lambda & z/r \\ 0 & 0 & 0 & \lambda
\end{array} \right), \; z= \pm 1, \; |r| > 1, \label{lemma3.2}
\end{equation} \begin{equation} H=D_{4}. \label{lemma3.3} \end{equation}
If $N$ has one real eigenvalue $\lambda$, $dim \:S_{0}=1$, $N_{1}$ is
decomposable, and $n=5$, then the pair $\{N,H\}$ is unitarily similar
to the canonical pair \{(\ref{lemma4.1}),(\ref{lemma4.2})\}:
\begin{equation} N=\left( \begin{array}{ccccc} \lambda & 1 & 0 &
\frac{1}{2}r^{2} & 0 \\ 0 & \lambda & 0 & z & 0 \\ 0 & 0 & \lambda & 0
& r \\ 0 & 0 & 0 & \lambda & 1 \\ 0 & 0 & 0 & 0 & \lambda \end{array}
\right), \; z= \pm 1, \; r>0, \label{lemma4.1} \end{equation}
\begin{equation} H=D_{5}.  \label{lemma4.2} \end{equation} If $N$ has one 
real eigenvalue $\lambda$, $dim \:S_{0}=1$, $N_{1}$ is decomposable,
and $n=6$, then the pair $\{N,H\}$ is unitarily similar to one and
only one of the canonical pairs \{(\ref{lemma5.1}),(\ref{lemma5.3})\},
\{(\ref{lemma5.2}),(\ref{lemma5.3})\}:
\begin{equation} N=\left( \begin{array}{cccccc} \lambda & 1   &   0 & 0  
    & 0 & 0 \\ 0 & \lambda & 1 & 0 & 0 & -r^{2}/2 \\ 0 & 0 & \lambda &
1 & 0 & 0 \\ 0 & 0 & 0 & \lambda & 0 & 1 \\ 0 & 0 & 0 & 0 & \lambda &
r \\ 0 & 0 & 0 & 0 & 0 & \lambda \end{array} \right),\; r>0,
\label{lemma5.1} \end{equation} \begin{equation} N=\left
( \begin{array}{cccccc} \lambda & 1 & -2r_{1} & 0 & 0 & 0 \\ 0 &
\lambda & 1 & r_{1} & 0 & -2r_{1}^{2}+r_{2}^{2}/2 \\ 0 & 0 & \lambda &
-1 & 0 & 0 \\ 0 & 0 & 0 & \lambda & 0 & -1 \\ 0 & 0 & 0 & 0 & \lambda
& r_{2} \\ 0 & 0 & 0 & 0 & 0 & \lambda \end{array} \right), \;
r_{2}>0, \label{lemma5.2} \end{equation}
\begin{equation} H=\left( \begin{array}{cccc} 0 & 0 & 0 & 1 \\ 0 & D_{3} 
& 0 & 0 \\ 0 & 0 & 1 & 0 \\ 1 & 0 & 0 & 0 \end{array} \right).
\label{lemma5.3}
\end{equation} If $N$ has one real eigenvalue $\lambda$, $dim \:S_{0}=2$, 
and $n=4$, then the pair $\{N,H\}$ is unitarily similar to one and
only one of the canonical pairs \{(\ref{lemma6.1}),(\ref{lemma6.5})\},
\{(\ref{lemma6.2}),(\ref{lemma6.5})\},
\{(\ref{lemma6.3}),(\ref{lemma6.5})\},
\{(\ref{lemma6.4}),(\ref{lemma6.5})\}: \begin{equation} N=\left (
\begin{array}{cccc} \lambda & 0 & \cos \alpha & \sin \alpha \\ 0 &
\lambda & - \sin \alpha & \cos \alpha \\ 0 & 0 & \lambda & 0 \\ 0 & 0 &
0 & \lambda \end{array} \right), \; 0< \alpha< \pi, \label{lemma6.1}
\end{equation} \begin{equation} N=\left( \begin{array}{cccc} \lambda &
0 & 0 & 1 \\ 0 & \lambda & r & 0 \\ 0 & 0 & \lambda & 0 \\ 0 & 0 & 0 &
\lambda \end{array} \right), \; |r|>1, \label{lemma6.2} \end{equation}
\begin{equation} N=\left( \begin{array}{cccc} \lambda & 0 &
\frac{1}{2} z & z \\ 0 & \lambda & -z & 0 \\ 0 & 0 & \lambda & 0 \\ 0
& 0 & 0 & \lambda \end{array} \right), \; z= \pm 1, \label{lemma6.3}
\end{equation} \begin{equation} N=\left( \begin{array}{cccc} \lambda &
0 & 0 & 0 \\ 0 & \lambda & 1 & 0 \\ 0 & 0 & \lambda & 0 \\ 0 & 0 & 0 &
\lambda \end{array} \right), \label{lemma6.4} \end{equation}
\begin{equation} H=\left( \begin{array}{cc} 0 & I_{2} \\ I_{2} & 0
\end{array} \right). \label{lemma6.5}
\end{equation} If $N$ has one real eigenvalue $\lambda$, $dim \:S_{0}=2$, 
and $n=5$, then the pair $\{N,H\}$ is unitarily similar to one and
only one of the canonical pairs \{(\ref{lemma7.1}),(\ref{lemma7.3})\},
\{(\ref{lemma7.2}),(\ref{lemma7.3})\}: \begin{equation} N=\left (
\begin{array}{ccccc} \lambda & 0 & 1 & 0 & 0 \\ 0 & \lambda & 0 & 1 &
0 \\ 0 & 0 & \lambda & z & 0 \\ 0 & 0 & 0 & \lambda & 0 \\ 0 & 0 & 0 &
0 & \lambda \end{array} \right), \; z= \pm 1, \label{lemma7.1}
\end{equation} \begin{equation} N=\left( \begin{array}{ccccc} \lambda
& 0 & 1 & 0 & 0 \\ 0 & \lambda & 0 & r & z \\ 0 & 0 & \lambda & 1 & 0
\\ 0 & 0 & 0 & \lambda & 0 \\ 0 & 0 & 0 & 0 & \lambda \end{array}
\right), \; z=\pm 1, \; r >0, \label{lemma7.2} \end{equation}
\begin{equation} H=\left( \begin{array}{ccc} 0 & 0 & I_{2} \\ 0 &
I_{1} & 0 \\ I_{2} & 0 & 0 \end{array} \right).  \label{lemma7.3}
\end{equation} If $N$ has one real eigenvalue $\lambda$, $dim
\:S_{0}=2$, and $n=6$, then the pair $\{N,H\}$ is unitarily similar to
one and only one of the canonical pairs
\{(\ref{lemma8.1}),(\ref{lemma8.3})\},
\{(\ref{lemma8.2}),(\ref{lemma8.3})\}: \begin{equation} N=\left
( \begin{array}{cccccc} \lambda & 0 & 1 & 0 & 0 & 0 \\ 0 & \lambda & 0
& 1 & r & 0 \\ 0 & 0 & \lambda & 0 & 1 & 0 \\ 0 & 0 & 0 & \lambda & 0
& 1 \\ 0 & 0 & 0 & 0 & \lambda & 0 \\ 0 & 0 & 0 & 0 & 0 & \lambda
\end{array} \right), \;\; r>0, \label{lemma8.1} \end{equation}
\begin{equation} N=\left( \begin{array}{cccccc} \lambda & 0 & 1 & 0 & 0 & 0 \\
 0 & \lambda & 0 & 1 & r & 0 \\ 0 & 0 & \lambda & 0 & \cos \alpha & \sin
\alpha \\ 0 & 0 & 0 & \lambda & -\sin \alpha & \cos \alpha \\ 0 & 0 & 0
& 0 & \lambda & 0 \\ 0 & 0 & 0 & 0 & 0 & \lambda \end{array} \right),
\;\; 0< \alpha< \pi, \label{lemma8.2} \end{equation} \begin{equation}
H=\left( \begin{array}{ccc} 0 & 0 & I_{2} \\ 0 & I_{2} & 0 \\ I_{2} &
0 & 0 \end{array} \right).  \label{lemma8.3}
\end{equation} If $N$ has one real eigenvalue $\lambda$, $dim \:S_{0}=2$, and 
$n=7$, then the pair $\{N,H\}$ is unitarily similar to the canonical
pair \{(\ref{lemma9.1}),(\ref{lemma9.2})\}:
\begin{equation} N=\left( \begin{array}{ccccccc} \lambda & 0 & 1 & 0 & 0 &
 0 & 0 \\ 0 & \lambda & 0 & 1 & 0 & 0 & 0 \\ 0 & 0 & \lambda & 0 & 0 &
\cos \alpha & -\sin \alpha \cos \beta \\ 0 & 0 & 0 & \lambda & 0 & \sin
\alpha & \cos \alpha \cos \beta \\ 0 & 0 & 0 & 0 & \lambda & 0 & \sin
\beta \\ 0 & 0 & 0 & 0 & 0 & \lambda & 0 \\ 0 & 0 & 0 & 0 & 0 & 0 &
\lambda \end{array} \right), \; 0< \alpha, \beta < \pi,
\label{lemma9.1} \end{equation} \begin{equation} H=\left(
\begin{array}{ccc} 0 & 0 & I_{2} \\ 0 & I_{3} & 0 \\ I_{2} & 0 & 0
\end{array} \right).  \label{lemma9.2} \end{equation} If $N$ has one
real eigenvalue $\lambda$, $dim \:S_{0}=2$, and $n=8$, then the pair
$\{N,H\}$ is unitarily similar to one and only one of the canonical
pairs \{(\ref{lemma10.1}),(\ref{lemma10.3})\},
\{(\ref{lemma10.2}),(\ref{lemma10.3})\}: $$N =\left(
\begin{array}{cccccccc} \lambda & 0 & 1 & 0 & 0 & 0 & 0 & 0 \\ 0 &
\lambda & 0 & 1 & 0 & 0 & 0 & 0 \\ 0 & 0 & \lambda & 0 & 0 & 0 & \cos
\alpha \sin \beta & \sin \alpha \sin \beta \\ 0 & 0 & 0 & \lambda & 0 & 0
& -\sin \alpha \sin \beta & \cos \alpha \sin \beta \\ 0 & 0 & 0 & 0 &
\lambda & 0 & \cos \beta & 0 \\ 0 & 0 & 0 & 0 & 0 & \lambda & 0 & \cos
\beta \\ 0 & 0 & 0 & 0 & 0 & 0 & \lambda & 0 \\ 0 & 0 & 0 & 0 & 0 & 0
& 0 & \lambda \end{array} \right), $$
\begin{equation} 0< \alpha<\pi, \;\; 0< \beta < \pi/2, \label{lemma10.1} 
\end{equation} $$ N =\left( \begin{array}{cccccccc}
\lambda & 0 & 1 & 0 & 0 & 0 & 0 & 0 \\ 0 & \lambda & 0 & 1 & 0 & 0 & 0
& 0 \\ 0 & 0 & \lambda & 0 & 0 & 0 & \cos \alpha \sin \beta & \sin \alpha
\sin \gamma \\ 0 & 0 & 0 & \lambda & 0 & 0 & -\sin \alpha \sin \beta &
\cos \alpha \sin \gamma \\ 0 & 0 & 0 & 0 & \lambda & 0 & \cos \beta & 0
\\ 0 & 0 & 0 & 0 & 0 & \lambda & 0 & \cos \gamma \\ 0 & 0 & 0 & 0 & 0 &
0 & \lambda & 0 \\ 0 & 0 & 0 & 0 & 0 & 0 & 0 & \lambda \end{array}
\right), $$ \begin{equation} 0< \alpha<\pi, \;\; 0 \leq \gamma< \beta
< \pi/2, \label{lemma10.2} \end{equation}
\begin{equation} H=\left( \begin{array}{ccc} 0 & 0 & I_{2} \\ 0 & I_{4} & 0 \\

I_{2} & 0 & 0 \end{array} \right).  \label{lemma10.3} \end{equation}
If $N$ has 2 distinct real eigenvalues $\lambda_{1}$, $\lambda_{2}$,
then the pair $\{N,H\}$ is unitarily similar to the canonical pair
\{(\ref{lemma11.1}),(\ref{lemma11.2})\}: \begin{equation} N=\left
( \begin{array}{cccc} \lambda_{1} & 1 & 0 & 0 \\ 0 & \lambda_{1} & 0 &
0 \\ 0 & 0 & \lambda_{2} & 0 \\ 0 & 0 & r & \lambda_{2}
\end{array} \right),\;\; for \; r \neq 0 \;\; \lambda_{1}< \lambda_{2},
 \label{lemma11.1} \end{equation} \begin{equation} H=\left(
\begin{array}{cc} 0 & I_{2} \\ I_{2} & 0 \end{array}
\right). \label{lemma11.2}
\end{equation} If $N$ has 3 eigenvalues: $\lambda \in R$, $\alpha \pm i 
\beta$ ($\alpha, \beta \in \Re$, $\beta > 0$), then the pair $\{N,H\}$
is unitarily similar to the canonical pair
\{(\ref{lemma12.1}),(\ref{lemma12.2})\}:
\begin{equation} N=\left( \begin{array}{cccc} \alpha & \beta & 0 &  0  
  \\ -\beta & \alpha & 0 &  0    \\
0 & 0 & \lambda & 0 \\ 0 & 0 & 0 & \lambda \end{array} \right),
\label{lemma12.1}
\end{equation} \begin{equation} H=\left( \begin{array}{cc} 0 & I_{2} \\
I_{2} & 0 \end{array} \right). \label{lemma12.2} \end{equation} If $N$
has 4 eigenvalues: $\alpha_{1} \pm i \beta_{1}$, $\alpha_{2} \pm i
\beta_{2}$, ($\alpha_{1}, \beta_{1}, \alpha_{2}, \beta_{2} \in \Re$,
$0< \beta_{1} \leq \beta_{2}$, $\alpha_{1} < \alpha_{2}$ if
$\beta_{1}=\beta_{2}$), then the pair $\{N,H\}$ is unitarily similar
to the canonical pair \{(\ref{lemma13.1}),(\ref{lemma13.2})\}:
\begin{equation} N=\left( \begin{array}{cccc} \alpha_{1} & \beta_{1} & 0 & 
 0 \\ -\beta_{1} & \alpha_{1} & 0 & 0 \\ 0 & 0 & \alpha_{2} & z
\beta_{2} \\ 0 & 0 & -z \beta_{2} & \alpha_{2}
\end{array} \right), \; z= \pm 1, \label{lemma13.1}
\end{equation} \begin{equation} H=\left( \begin{array}{cc} 0 & I_{2} \\
I_{2} & 0 \end{array} \right). \label{lemma13.2} \end{equation} If $N$
has 2 eigenvalues $\alpha \pm i \beta$ ($\alpha, \beta \in \Re$,
$\beta>0$), and $n=4$, then the pair $\{N,H\}$ is unitarily similar to
one and only one of the canonical pairs
\{(\ref{lemma14.1}),(\ref{lemma14.3})\},
\{(\ref{lemma14.2}),(\ref{lemma14.3})\}: \begin{equation} N=\left
( \begin{array}{cccc} \alpha & \beta & \cos \gamma & \sin \gamma \\
-\beta & \alpha & -\sin \gamma & \cos \gamma \\ 0 & 0 & \alpha & \beta
\\ 0 & 0 & -\beta & \alpha \end{array} \right), \; 0 \leq \gamma < 2
\pi, \label{lemma14.1} \end{equation} \begin{equation} N=\left
( \begin{array}{cccc} \alpha & \beta & 0 & 1 \\ -\beta & \alpha & 1 &
0 \\ 0 & 0 & \alpha & -\beta \\ 0 & 0 & \beta & \alpha \end{array}
\right), \label{lemma14.2} \end{equation} \begin{equation} H=\left
( \begin{array}{cc} 0 & I_{2} \\ I_{2} & 0 \end{array}
\right). \label{lemma14.3}
\end{equation} If $N$ has 2 eigenvalues $\alpha \pm i \beta$ ($\alpha, 
\beta \in \Re$, $\beta > 0$), and $n=6$, then the pair $\{N,H\}$ is
unitarily similar to one and only one of the canonical pairs
\{(\ref{lemma15.1}),(\ref{lemma15.3})\},
\{(\ref{lemma15.2}),(\ref{lemma15.3})\}: $$ N= \left (
\begin{array}{cccccc} \alpha & \beta & 0 & 0 & 0 & r \\ -\beta &
\alpha & 0 & 1 & (\cos \gamma + 1)/4 \beta - r & \sin \gamma/4 \beta \\
0 & 0 & \alpha & \beta & \frac{1}{2}(\cos \gamma + 1) & \frac{1}{2}\sin
\gamma \\ 0 & 0 & -\beta & \alpha & -\frac{1}{2}\sin \gamma &
\frac{1}{2}(\cos \gamma -1) \\ 0 & 0 & 0 & 0 & \alpha & \beta \\ 0 & 0
& 0 & 0 & -\beta & \alpha \end{array} \right), $$ \begin{equation} 0
\leq \gamma < 2\pi, \; \gamma \neq \pi, \label{lemma15.1}
\end{equation}
\begin{equation} N= \left( \begin{array}{cccccc} \alpha & \beta & 0 & 0 
& r & 0 \\ -\beta & \alpha & 0 & 1 & 0 & r \\ 0 & 0 & \alpha & \beta &
0 & 0 \\ 0 & 0 & -\beta & \alpha & 0 & -1 \\ 0 & 0 & 0 & 0 & \alpha &
\beta \\ 0 & 0 & 0 & 0 & -\beta & \alpha \end{array} \right),
\label{lemma15.2} \end{equation} \begin{equation} H=\left
( \begin{array}{ccc} 0 & 0 & I_{2} \\ 0 & I_{2} & 0 \\ I_{2} & 0 & 0
\end{array} \right). \label{lemma15.3}
\end{equation} If $N$ has 2 eigenvalues $\alpha \pm i \beta$ ($\alpha, 
\beta \in \Re$, $\beta > 0$), and $n=8$, then the pair $\{N,H\}$ is
unitarily similar to the canonical pair
\{(\ref{lemma16.1}),(\ref{lemma16.2})\}:
$$ N= \left( \begin{array}{cccccccc} \alpha & \beta & 0 & 0 & 0 & 0 &
0 & 0 \\ -\beta & \alpha & 0 & 1 & 0 & 0 & \sin^{2}\gamma/2\beta & \sin
\gamma \cos \gamma \cos \delta /2 \beta \\ 0 & 0 & \alpha & \beta & 0 &
0 & \sin^{2} \gamma & \sin \gamma \cos \gamma \cos \delta \\ 0 & 0 &
-\beta & \alpha & 0 & 0 & -\sin \gamma \cos \gamma \cos \delta & -\cos^{2}
\gamma \\ 0 & 0 & 0 & 0 & \alpha & \beta & \sin \gamma \cos \gamma \sin
\delta & 0 \\ 0 & 0 & 0 & 0 & -\beta & \alpha & 0 & \sin \gamma \cos
\gamma \sin \delta \\ 0 & 0 & 0 & 0 & 0 & 0 & \alpha & \beta \\ 0 & 0 &
0 & 0 & 0 & 0 & -\beta & \alpha \end{array} \right) $$
\begin{equation} 0<\gamma<\pi/2, \; 0<\delta <\pi, \label{lemma16.1}
\end{equation} \begin{equation} H=\left( \begin{array}{ccc} 0 & 0 & I_{2} \\
0 & I_{4} & 0 \\ I_{2} & 0 & 0 \end{array} \right). \label{lemma16.2}
\end{equation} Here all parameters are $H$-unitary invariants, i.e.,
the same canonical forms are $H$-unitarily similar to each other iff
the values of all parameters are equal.  \end{theorem}

The proof of the theorem is presented in what follows.

\subsection{One Real Eigenvalue of $N$}

The case when $N$ has only one real eigenvalue $\lambda$ can be considered
as in~\cite{2}. Namely, if $dim\:S_{0}=1$, then there exists two alternatives:
$N_{1}$ is indecomposable or decomposable, this property being independent
of the choice of the decomposition $R^{n}=S_{0}\dot{+}S\dot{+}S_{1}$
because the indecomposability or decomposability of $N_{1}$ does not change under
unitary similarity of the pair $\{N_{1},H_{1}\}$. In the former case
one can show that $n \leq 5$ and obtain the canonical forms
\{(\ref{lemma1.1}),(\ref{lemma1.2})\} - \{(\ref{lemma2.2}),(\ref{lemma2.3})\},
in the latter one can show that $n \leq 6$ and obtain the canonical forms
\{(\ref{lemma3.1}),(\ref{lemma3.3})\} - \{(\ref{lemma5.2}),(\ref{lemma5.3})\}
in just the same way as it was done in \cite{2}. If the subspace $S_{0}$
is two-dimensional, the operator $N$ can also be considered as in~\cite{2}
except for the
case $n=4$ because one of the corresponding canonical forms in~\cite{2} is
essentially complex. Thus, for the case when $N$ has one real eigenvalue
$\lambda$ we will consider only the alternative $dim\:S_{0}=2$, $n=4$ and omit
the rest.

\subsubsection{$dim\:S_{0}=2$, $n=4$}
In this case $R^{4}=S_{0}\dot{+}S_{1}$. Therefore,
$$N - \lambda I=\left( \begin{array}{cc}
0 & N_{2} \\
0 & 0 \end{array}
\right)=
\left( \begin{array}{cccc}
0 & 0 & a & b \\
0 & 0 & c & d \\
0 & 0 & 0 & 0 \\
0 & 0 & 0 & 0 \end{array}
\right), \;\;
H=\left( \begin{array}{cc}
0 & I_{2} \\
I_{2} & 0 \end{array}
\right),
$$
and the submatrix $N_{2}$ is not restricted by the condition of the
$H$-normality of $N$.

(a) $det N_{2} \neq 0.\/$ Suppose an $H$-unitary transformation
$$T=\left( \begin{array}{cc}
T_{1} & T_{2} \\
T_{3} & T_{4} \end{array}
\right)$$
reduces $N- \lambda I$ to the form $\tilde{N}-\lambda I$:
$$ N-\lambda I=\left( \begin{array}{cc}
0 & N_{2} \\
0 & 0 \end{array}
\right),\;\;\;
\tilde{N}-\lambda I=\left( \begin{array}{cc}
0 & \widetilde{N_{2}} \\
0 & 0 \end{array}
\right). $$
Then conditions (\ref{eqn1}) - (\ref{eqn3}) below are necessarily satisfied:
\begin{eqnarray}
N_{2}T_{3} & = & 0,  \label{eqn1} \\
N_{2}T_{4} & = & T_{1} \widetilde{N_{2}}, \label{eqn2} \\
0 & = & T_{3} \widetilde{N_{2}}. \label{eqn3}
\end{eqnarray}
Since $N_{2}$ is nondegenerate, (\ref{eqn1}) is satisfied only if
$T_{3}=0$. The operator $T$ is $H$-unitary iff
\begin{eqnarray}
T_{1}T_{4}^{*} & = & I,  \label{hh1} \\
T_{1}T_{2}^{*}+T_{2}T_{1}^{*} & = & 0. \label{hh2}
\end{eqnarray}
It follows from system (\ref{hh1}) - (\ref{hh2}) that without loss of
generality we can consider only quasidiagonal transformations
$T=T_{1} \oplus T_{1}^{*-1}$ because $T_{2}$ does not appear in equations
(\ref{eqn1}) - (\ref{eqn3}).

Thus, the only condition
\begin{equation}
N_{2}=T_{1}\widetilde{N_{2}}T_{1}^{*}  \label{meq}
\end{equation}
should be satisfied, i.e., it is necessary to find out what form a
nondegenerate $2 \times 2$-matrix $N_{2}$ can be reduced to under congruence.

Consider the matrix $N_{2}'=N_{2}N_{2}^{*-1}$. Its spectral characteristics
are invariant because $N_{2}'=T_{1}\widetilde{N_{2}}'T_{1}^{-1}$. Since
$det\: N_{2}'=1$, $N_{2}'$ has either two complex conjugate eigenvalues
$\cos \alpha \pm i \sin \alpha$ or two real eigenvalues $r$, $1/r$ ($r \neq 0$).
In the former case $N_{2}'$ can be reduced to the form
\begin{equation}
N_{2}'=\left( \begin{array}{cc}
\cos \alpha & \sin \alpha \\
-\sin \alpha & \cos \alpha \end{array}
\right), \;\; 0 < \alpha< \pi, \label{first4}
\end{equation}
in the latter to the Jordan normal form.

If $N_{2}'$ has form (\ref{first4}), then
$$N_{2}=\left( \begin{array}{cc}
t \: \sin \alpha/(1-\cos \alpha) & t \\
-t & t \: \sin \alpha/(1-\cos \alpha) \end{array}
\right), \; t \neq 0.
$$
As $det\: N_{2}=2t^{2}/(1-\cos \alpha) >0$, one can take
$T_{1}=\sqrt{det \: N_{2}} I $ and obtain
$$N_{2}=\left( \begin{array}{cc}
\pm \cos \frac{\alpha}{2} & \pm \sin \frac{\alpha}{2} \\
\mp \sin \frac{\alpha}{2} & \pm \cos \frac{\alpha}{2} \end{array}
\right), \;\; 0 < \alpha < \pi.
$$
Since the transformation $T_{1}=D_{2}$
replaces $\sin \frac{\alpha}{2}$ by $- \sin \frac{\alpha}{2}$,
we can write
\begin{equation}
N_{2}=\left( \begin{array}{cc}
\cos \alpha & \sin \alpha \\
-\sin \alpha & \cos \alpha \end{array}
\right), \;\; 0 < \alpha < \pi \label{form41}
\end{equation}
(note that two last formulas for $N_{2}$ are not equivalent because
(\ref{form41}) includes the extra value $\alpha=\pi/2$ corresponding
to the case $N_{2}'=-I$).

Now we must prove the invariance of the parameter $\alpha$. To this end
suppose that a nondegenerate matrix $T_{1}$ satisfies (\ref{meq}), where
$N_{2}$ has form (\ref{form41}) and
$$\widetilde{N_{2}}=\left( \begin{array}{cc}
\cos \tilde{\alpha} & \sin \tilde{\alpha} \\
-\sin \tilde{\alpha} & \cos \tilde{\alpha} \end{array}
\right), \;\; 0 < \tilde{\alpha} < \pi.
$$
As $N_{2}+N_{2}^{*}=T_{1}(\widetilde{N_{2}}+\widetilde{N_{2}}^{*})T_{1}^{*}$
and $N_{2}-N_{2}^{*}=T_{1}(\widetilde{N_{2}}-\widetilde{N_{2}}^{*})T_{1}^{*}$,
we have 
$$ 4 \cos^{2} \alpha=det (N_{2}+N_{2}^{*})=(det T_{1})^{2} \: det 
(\widetilde{N_{2}}+\widetilde{N_{2}^{*}})=(det T_{1})^{2}4\cos^{2}
\tilde{\alpha} $$
and $$ 4 \sin^{2} \alpha=det (N_{2}-N_{2}^{*})=(det T_{1})^{2} \: det (\widetilde{N_{2}}-\widetilde{N_{2}^{*}})=(det T_{1})^{2}4\sin^{2}\tilde{\alpha}.$$
Therefore, $|det T_{1}|=1$, $\cos \alpha = \pm \cos \tilde{\alpha}$, $\sin \alpha=\sin \tilde{\alpha}$.
Now we write the condition
$N_{2}+N_{2}^{*}=T_{1}(\widetilde{N_{2}}+\widetilde{N_{2}}^{*})T_{1}^{*}$
in detail:
$$\left( \begin{array}{cc}
\cos \alpha & 0 \\
0 & \cos \alpha \end{array}
\right)=\left( \begin{array}{cc}
(t_{11}^{2}+t_{12}^{2})\cos \tilde{\alpha} & (t_{11}t_{21}+t_{12}t_{22})\cos \tilde{\alpha} \\
(t_{11}t_{21}+t_{12}t_{22})\cos \tilde{\alpha} & (t_{21}^{2}+t_{22}^{2})\cos \tilde{\alpha} \end{array}
\right). $$
Since $|\cos \alpha|=|\cos \tilde{\alpha}|$, we have $t_{11}^{2}+t_{12}^{2}=1$, hence
$\cos \alpha=\cos \tilde{\alpha}$. Thus, $\alpha=\tilde{\alpha}$, Q.E.D.

If $N_{2}'$ has distinct real eigenvalues $r$ and $1/r$, i.e.,
$r \neq \pm 1$, then it can be reduced to the diagonal form
$N_{2}'=1/r \oplus r$, $\: |r|>1$.
Consequently,
$$ N_{2}=\left( \begin{array}{cc}
0 & t \\
rt & 0 \end{array}
\right), \; t \neq 0. $$
Taking $T_{1}=1 \oplus t$,
we reduce $N_{2}$ to the form
\begin{equation}
N_{2}=\left( \begin{array}{cc}
0 & 1 \\
r & 0 \end{array}
\right), \; |r|>1. \label{form42}
\end{equation}
It is clear that $r$ is an invariant.

Finally, we consider the case when $N_{2}'$ has the eigenvalues $\pm 1$.
If $N_{2}'=I$, the matrix $N_{2}$ is selfadjoint, hence, it can be reduced to
the diagonal form. Therefore, the nondegenerate subspace
$V=span \{ v_{1}, v_{3} \}$ is invariant both for $N$ and for $N^{[*]}$, i.e.,
the operator $N$ is decomposable. It can easily be checked that $N_{2}'$
is not equivalent to the form
$$ N_{2}'=\left( \begin{array}{cc}
1 & 1 \\
0 & 1 \end{array} \right) $$
because then $N_{2}$ turns out to be degenerate, which is impossible.
If $N_{2}'=-I$, $N_{2}$ can be reduced to the above-mentioned form
$$N_{2}=\left( \begin{array}{cc}
0 & 1 \\
-1 & 0 \end{array}
\right).$$
The last case to be considered is the case when the Jordan normal
form of $N_{2}'$ is
$$N_{2}'=\left( \begin{array}{cc}
-1 & 1 \\
0 & -1 \end{array}
\right).$$
Then
$$N_{2}=\left( \begin{array}{cc}
\frac{1}{2}t & t \\
-t & 0 \end{array}
\right) , \; t \neq 0. $$
Taking $T_{1}=\sqrt{|t|}I$, we achieve
\begin{equation}
N_{2}=\left( \begin{array}{cc}
\frac{1}{2}z & z \\
-z & 0 \end{array}
\right), \; z=\pm 1.  \label{form43}
\end{equation}
Here $z$ is an invariant. Indeed, suppose that some matrix $T_{1}$
satisfies condition (\ref{meq}), where
$$\widetilde{N_{2}}=\left( \begin{array}{cc}
\frac{1}{2}\tilde{z} & \tilde{z} \\
-\tilde{z} & 0 \end{array}
\right), \; \tilde{z}=\pm 1. $$
Then $\frac{1}{2}z=\frac{1}{2}t_{11}^{2}\tilde{z}$, hence $z=\tilde{z}$.

As a result, we have obtained three forms (\ref{form41}), (\ref{form42}),
(\ref{form43}). Now it is necessary to find out whether the operator $N$ is
indecomposable in the three cases. The indecomposability of $N$ means that
$(aN_{2}+bN_{2}^{*})x=0$ only if $(x,N_{2}x)=0$ ($a^{2}+b^{2} \neq 0$).
If $N_{2}'=N_{2}N_{2}^{*-1}$ has no real eigenvalues, the equation
$(aN_{2}+bN_{2}^{*})x=0$ has no solutions, i.e., $N$ is indecomposable if
$N_{2}$ has form (\ref{form41}) with $\alpha \neq \pi/2$. If an eigenvalue
$\lambda$ of $N_{2}'$ is not equal to 1, then $(x,N_{2}x)=0$ because
$(x,N_{2}x)=(x,\lambda N_{2}^{*}x)=\lambda (x, N_{2}^{*}x)=\lambda (x,N_{2}x)$.
Thus, if $N_{2}$ has form (\ref{form42}), (\ref{form43}), or (\ref{form41}) with
$\alpha = \pi/2$, then $N$ is also indecomposable.

(b) $det N_{2}=0.\/$ Since $N$ with $N_{2}=0$ is decomposable, it suffices to
consider the remaining case $rg\:N_{2}=1$:
$$N_{2}=\left( \begin{array}{cc}
ka & kb \\
la & lb \end{array}
\right), \; a^{2}+b^{2} \neq 0, \; k^{2}+l^{2} \neq 0.$$
It is readily seen that $S_{0} \cap S_{1} \neq \{0\}$ if $la=kb$, therefore, we can assume
that this condition is not satisfied. Taking $T=T_{1}\oplus T_{1}^{*-1}$,
where
$$T_{1}=\left( \begin{array}{cc}
a & k  \\
b & l \end{array}
\right), $$
we obtain one more canonical form:
$$N - \lambda I=\left( \begin{array}{cccc}
0 & 0 & 0 & 0 \\
0 & 0 & 1 & 0 \\
0 & 0 & 0 & 0 \\
0 & 0 & 0 & 0 \end{array}
\right)$$
(it can easily be checked that this form is indecomposable).

As a result, we have proved that

{\em{If an indecomposable $H$-normal operator $N$
($N: \: C^{4} \rightarrow C^{4}$) has the only eigenvalue $\lambda \in \Re$,
and $dim \:S_{0}=2$, then the pair $\{N,H\}$ is unitarily similar to one and
only one of the canonical pairs
\mbox{\{(\ref{lemma6.1}),(\ref{lemma6.5})\}},
\mbox{\{(\ref{lemma6.2}),(\ref{lemma6.5})\}},
\mbox{\{(\ref{lemma6.3}),(\ref{lemma6.5})\}},
\mbox{\{(\ref{lemma6.4}),(\ref{lemma6.5})\}}.}}

\subsection{Two Real Distinct Eigenvalues of $N$}
Since the canonical pair \{(\ref{lemma11.1}),(\ref{lemma11.2})\} is obtained
in the same way as in~\cite{2}, we will not repeat the proof of the
following fact:

{\em{If an indecomposable $H$-normal operator acts in
a space $R^{n}$ of rank 2 and has 2 distinct real eigenvalues: $\lambda_{1}$
and $\lambda_{2}$, then $n=4$ and
the pair $\{N,H\}$ is unitarily similar to
the canonical pair \{(\ref{lemma11.1}),(\ref{lemma11.2})\}.}}

\subsection{Three Eigenvalues of $N$: One Real and Two Complex Conjugate}
Suppose an indecomposable $H$-normal operator $N$ has a real eigenvalue
$\lambda$ and two complex eigenvalues $\alpha \pm i \beta$ ($\beta > 0$).
According to Lemma~2.1, we have $R^{n}={\cal Q}_{1}\dot{+}{\cal Q}_{2}$,
$dim\,{\cal Q}_{1}=dim\,{\cal Q}_{2}=m$,
$[{\cal Q}_{1},{\cal Q}_{1}]=0$, $[{\cal Q}_{2},{\cal Q}_{2}]=0$,
$N{\cal Q}_{1} \subseteq {\cal Q}_{1}$,
$N{\cal Q}_{2} \subseteq {\cal Q}_{2}$,
$N_{1}=N|_{{\cal Q}_{1}}$ has two eigenvalues $\alpha \pm i \beta$,
$N_{2}=N|_{{\cal Q}_{2}}$ one eigenvalue $\lambda$.
Since $min\{v_{-},v_{+}\}=2$, $n=v_{-}+v_{+} \geq 4$.
On the other hand, the subspaces ${\cal Q}_{1}$ and ${\cal Q}_{2}$
are neutral so that $n=2m \leq 4$. Thus, $n=4$. As $H$ is nondegenerate,
for any basis in ${\cal Q}_{1}$ there exists a basis in ${\cal Q}_{2}$
such that
$$H=\left( \begin{array}{cc}
0 & I  \\
I & 0
\end{array}
\right). $$
Take a basis in ${\cal Q}_{1}$ such that
\begin{equation}
N_{1}=\left( \begin{array}{cc}
\alpha & \beta  \\
-\beta & \alpha	\end{array}
\right).  \label{n1}
\end{equation}
Then with respect to the decomposition $R^{n}={\cal Q}_{1}\dot{+}{\cal Q}_{2}$
we have
\begin{equation}
N=\left( \begin{array}{cc}
N_{1} & 0  \\
0 & N_{2}
\end{array}
\right), \;\;
H=\left( \begin{array}{cc}
0 & I  \\
I & 0
\end{array}
\right). \label{n_and_h}
\end{equation}
The condition of the $H$-normality of $N$ is
\begin{equation}
N_{1}N_{2}^{*}=N_{2}^{*}N_{1}. \label{h_norm2}
\end{equation}
The only matrix commuting with (\ref{n1}) and having one eigenvalue
$\lambda$ is $\lambda I$. Thus,
$$N=\left( \begin{array}{cc}
\alpha & \beta   \\
-\beta & \alpha  \end{array} \right) \oplus
\left( \begin{array}{cc}
\lambda & 0 \\
0 & \lambda \end{array}
\right). $$
It can easily be checked that $N$ is indecomposable. Indeed, suppose a subspace
$V$ is invariant for $N$ and $N^{[*]}$. Since
$min \{ dim \:V, dim \: V^{[\perp]}\} \leq 2$,
we can assume that $dim \: V \leq 2$. If $V$ were of dimension 1, then there
would exist a vector $v \in V$ such that $Nv=\lambda v$, $N^{[*]}v=\lambda v$.
But all eigenvectors of $N$ corresponding to the eigenvalue $\lambda$ are not
eigenvectors of $N^{[*]}$. Thus, $dim\:V \neq 1$. Suppose $dim \:V=2$. Then
$N|_{V}$ has either the only eigenvalue $\lambda$ or two eigenvalues
$\alpha \pm i \beta$. In the former case $V={\cal Q}_{2}$, in the latter
$V={\cal Q}_{1}$. In the both cases $V$ is degenerate, therefore, $N$ is
indecomposable.

Thus, we have proved that

{\em{If an indecomposable $H$-normal operator acts in
a space $R^{n}$ of rank 2 and has 3 eigenvalues: $\lambda \in R$,
$\alpha \pm i \beta$ ($\alpha, \beta \in \Re$, $\beta > 0$), then $n=4$ and
the pair $\{N,H\}$ is unitarily similar to the
canonical pair \{(\ref{lemma12.1}),(\ref{lemma12.2})\}.}}

\subsection{Two Distinct Pairs of Complex Conjugate Eigenvalues of $N$}
Suppose $N$ has four eigevalues $\alpha_{1} \pm i \beta_{1}$,
$\alpha_{2} \pm i \beta_{2}$ ($\beta_{1}, \; \beta_{2}>0$,
$(\alpha_{1}, \beta_{1}) \neq (\alpha_{2},\beta_{2})$). Let us fix the order
of these pairs: $\beta_{1} \leq \beta_{2}$, $\alpha_{1} < \alpha_{2}$ if
$\beta_{1}=\beta_{2}$. As in the previous section, one can show that $N$
and $H$ can be reduced to form (\ref{n_and_h}) with
$$N_{1}=\left( \begin{array}{cc}
\alpha_{1} & \beta_{1}  \\
-\beta_{1} & \alpha_{1}	\end{array}
\right). $$
It follows from condition (\ref{h_norm2}) of the $H$-normality of $N$ that
$$N_{2}=\left( \begin{array}{cc}
\alpha_{2} &  z \beta_{2}  \\
-z \beta_{2} & \alpha_{2} \end{array}
\right), \; z=\pm 1. $$
Now we prove that the number $z$ is an $H$-unitary invariant. To this end
suppose that a matrix $T$ satisfies condition (\ref{similar}) $NT=T\tilde{N}$
and condition (\ref{unitar}) $TT^{[*]}=I$, where
$$N=N_{1} \oplus \left( \begin{array}{cc}
\alpha_{2} & z \beta_{2} \\
-z \beta_{2} & \alpha_{2} \end{array}
\right), \;\;
\tilde{N}=N_{1} \oplus \left( \begin{array}{cccc}
\alpha_{2} & \tilde{z} \beta_{2} \\
-\tilde{z} \beta_{2} & \alpha_{2} \end{array}
\right),\;\; |z|=|\tilde{z}|=1. $$
It follows from (\ref{similar}) that $T=T_{1}\oplus T_{2}$, where
$$ T_{1}=\left( \begin{array}{cc}
t_{11} & t_{12}  \\
-t_{12} & t_{11} \end{array}
\right). $$
It follows from (\ref{unitar}) that $T_{2}=T_{1}^{*-1}$, therefore,
$$T_{2}=\left( \begin{array}{cc}
t_{33} & t_{34}  \\
-t_{34} & t_{33} \end{array}
\right). $$
It is seen that under these conditions $\tilde{z}=z$.
The indecomposability of the form obtained can be checked as before.

Thus, we have proved that

{\em{If an indecomposable $H$-normal operator acts in a space $R^{n}$
of rank 2 and has 4 eigenvalues: $\alpha_{1} \pm i \beta_{1}$,
$\alpha_{2} \pm i \beta_{2}$, ($\alpha_{1}, \beta_{1}, \alpha_{2},
\beta_{2} \in \Re$, $0< \beta_{1} \leq \beta_{2}$,
$\alpha_{1} < \alpha_{2}$ if $\beta_{1}=\beta_{2}$),
then $n=4$ and the pair $\{N,H\}$ is unitarily similar to the
canonical pair \{(\ref{lemma13.1}),(\ref{lemma13.2})\}.}}

\subsection{Two Complex Conjugate Eigenvalues of $N$}
The two following propositions hold for any space with indefinite scalar
product. They are in a sense analogous to Propositions 1, 2
from~\cite{2}.

\begin{proposition}
Let an indecomposable $H$-normal operator $N$ acting in $R^{n}$
($n>2$) have two distinct eigenvalues $\lambda=\alpha+i\beta$,
$\overline{\lambda}=\alpha-i\beta$. Let
$$S_{0}'=\{z=x+iy \; (x,y \in R^{n}): \;  Nz=\lambda z, \; N^{[*]}z=\overline{\lambda}z \}, $$
$$S_{0}''=\{z=x+iy \; (x,y \in R^{n}): \;  Nz=\lambda z, \; N^{[*]}z=\lambda z \}, $$
$\{z_{j}\}_{1}^{p}$ ($\{z_{j}\}_{p+1}^{p+q}$) be a basis of $S_{0}'$ ($S_{0}''$),
and
$$S_{0}=\sum_{j=1}^{p+q}span\{x_{j},y_{j}\}.$$
Then there exists a decomposition of $R^{n}$ into a direct sum of
subspaces $S_{0}$, $S$, $S_{1}$ such that
\begin{equation}
N=\left( \begin{array}{ccc}
N' & * & *    \\
0 & N_{1} & * \\
0 & 0 & N''
\end{array} \right), \;\;
H=\left( \begin{array}{ccc}
0 & 0 & I    \\
0 & H_{1} &  0    \\
I & 0 & 0 \end{array} \right), \label{aug_dec1}
\end{equation}
where
$$ N': \; S_{0} \rightarrow S_{0}, \;\; N'=N_{1}' \oplus \ldots \oplus N_{p+q}', $$
\begin{equation}
N_{j}'=\left( \begin{array}{cc}
\alpha & \beta    \\
-\beta & \alpha \end{array} \right), \;\;  j=1, \ldots p+q,
\label{aug_dec2}
\end{equation}
$$ N'': \; S_{1} \rightarrow S_{1}, \;\; N''=N_{1}'' \oplus \ldots \oplus N_{p+q}'', $$
\begin{equation}
N_{j}''=N_{j}' \;\; if \;1 \leq j \leq p, \;\;
N_{j}''=N_{j}'^{*} \;\; if \; p < j \leq p+q, \label{aug_dec3}
\end{equation}
the internal operator $N_{1}$ is $H_{1}$-normal and the pair
$\{N_{1},H_{1}\}$ is determined up to unitarily similarity.
To go over from one decomposition $R^{n}=S_{0}\dot{+}S\dot{+}S_{1}$
to another by means of a transformation $T$ it is necessary that the
matrix $T$ be block triangular with respect to both decompositions.
\end{proposition}

{\bf{Proof:}} It is clear that the subspace $S_{0}$ is well defined,
i.e., that its definition does not depend on the choice of bases in
$S_{0}'$ and $S_{0}''$. Since $N$ and $N^{[*]}$ commute and have
two eigenvalues, at least one of the subspaces $S_{0}'$, $S_{0}''$
is nontrivial so that $p+q >0$. Show that the system
$\{x_{j}\}_{1}^{p+q} \cup \{y_{j}\}_{1}^{p+q}$ is a basis in
$S_{0}$. In fact, the assumption
$\sum_{j=1}^{p+q}(a_{j}x_{j}+b_{j}y_{j})=0$
($a_{j}, b_{j} \in \Re,$ $j=1, \ldots p+q$)
means that ${\cal R}e \sum_{j=1}^{p+q}(a_{j} - ib_{j})z_{j}=0$,
therefore, ${\cal R}e \{N \sum_{j=1}^{p+q}(a_{j} - ib_{j})z_{j}\}=0$.
But
${\cal R}e \{N \sum_{j=1}^{p+q}(a_{j} - ib_{j})z_{j}\}=\alpha{\cal R}e\sum_{j=1}^{p+q}(a_{j}-ib_{j})z_{j}-\beta{\cal I}m\sum_{j=1}^{p+q}(a_{j}-ib_{j})z_{j}$
so that ${\cal I}m \sum_{j=1}^{p+q}(a_{j} - ib_{j})z_{j}=0$.
Thus, $\sum_{j=1}^{p+q}(a_{j} - ib_{j})z_{j}=0$. Since the vectors
$z_{j}$ are linearily independent in $C^{n}$, $a_{j}=b_{j}=0$
($j=1, \ldots p+q$), i.e., the vectors
$\{x_{j}\}_{1}^{p+q} \cup \{y_{j}\}_{1}^{p+q}$ are linearily independent
in $R^{n}$.  Thus, the dimension of $S_{0}$ is equal to $2(p+q)$.

Now let us prove that for $N$ to be indecomposable it is necesssary that
$S_{0}$ be neutral. Indeed, we already know that if $z=x+iy$
($x,y \in R^{n}$) is an eigenvector of $N^{[*]}$ such that $Nz=\lambda z$,
then the subspace $span\{x,y\}$, which is invariant for $N$ and $N^{[*]}$,
is either nondegenerate or neutral (see Section~2.3). Since
$n>2$ and $N$ is indecomposable, it is necessarily neutral.
Further, if
$Nz_{1}=\lambda z_{1}$, $N^{[*]}z_{1}=\overline{\lambda}z_{1}$,
$Nz_{2}=\lambda z_{2}$, $N^{[*]}z_{2}=\lambda z_{2}$,
then it can be shown (as in Section~2.3) that
$[z_{1},z_{2}]=[z_{1},\overline{z_{2}}]=0$, hence
$[x_{1},x_{2}]=[x_{1},y_{2}]=[y_{1},x_{2}]=[y_{1},y_{2}]=0$.
If
$Nz_{1}=\lambda z_{1}$, $N^{[*]}z_{1}=\overline{\lambda}z_{1}$,
$Nz_{2}=\lambda z_{2}$, $N^{[*]}z_{2}=\overline{\lambda}z_{2}$,
then $[z_{1},\overline{z_{2}}]=0$, i.e.,
$[x_{1},x_{2}]=[y_{1},y_{2}]$ and $[x_{1},y_{2}]=-[y_{1},x_{2}]$.
If $a^{2}+b^{2} \neq 0$ ($a=[x_{1},x_{2}]$, $b=[x_{1},y_{2}]$),
the two-dimensional subspace
$span\{ax_{1}-by_{1}+x_{2}, bx_{1}+ay_{1}+y_{2}\}$, which is
invariant for $N$ and $N^{[*]}$,
will be nondegenerate, therefore, $N$ will be decomposable.
Thus, for $N$ to be indecomposable it is necessary to have
$a=b=0$. It can be checked in the similar way that
the conditions
$[x_{1},x_{2}]=[y_{1},y_{2}]=[x_{1},y_{2}]=[y_{1},x_{2}]=0$
are satisfied if
$Nz_{1}=\lambda z_{1}$, $N^{[*]}z_{1}={\lambda}z_{1}$,
$Nz_{2}=\lambda z_{2}$, $N^{[*]}z_{2}={\lambda}z_{2}$.
Thus, if $N$ is indecomposable, $S_{0}$ is neutral.

For any neutral subspace $S_{0}$ of a space with indefinite scalar
product there exists a subspace $S_{1}$ such that
$$   H|_{(S_{0}\dot{+}S_{1})}=\left( \begin{array}{cc}
0 & I    \\
I & 0  \end{array} \right). $$
Since $(S_{0} \dot{+} S_{1})$ is nondegenerate, the subspace
$S=(S_{0} \dot{+} S_{1})^{[\perp]}$ is nondegenerate too and
$R^{n}=S_{0} \dot{+}S \dot{+}S_{1}$.
It is clear that with respect to this decomposition
the matrices $N$ and $H$ have form (\ref{aug_dec1}), the submatrix
$N'$ has form (\ref{aug_dec2}) and $N''$ has from (\ref{aug_dec3}).
The last two statements of the proposition can be proved as in
Proposition~1 from~\cite{2}. The proof is completed.

\begin{proposition}
An $H$-normal operator such that $dim \:S_{0}=2$ is indecomposable.
\end{proposition}
{\bf{Proof:}} Assume the converse.
Suppose some nondegenerate subspace $V$ is invariant both for $N$ and for
$N^{[*]}$. Let us denote $V_{1}=V$, $V_{2}=V^{[\perp]}$, $N_{1}=N|_{V_{1}}$,
$N_{2}=N|_{V_{2}}$, $H_{1}=H|_{V_{1}}$, $H_{2}=H|_{V_{2}}$.
Since the operators $N_{i}$ ($i=1,2$) are $H_{i}$-normal,
both subspaces $S_{0}^{(i)} \subset V_{i}$ (defined as $S_{0}$) are nontrivial,
i.e., $dim \:S_{0}^{(i)} \geq 2$.
Since $S_{0}=S_{0}^{(1)}\dot{+}S_{0}^{(2)}$,
$dim\:S_{0}=dim\:S_{0}^{(1)}+dim\:S_{0}^{(2)} \geq 4$.
This contradicts the condition $dim\:S_{0}=2$. Thus, $N$ is indecomposable.

Now let us show that if $min\{v_{-},v_{+}\}=2$, then $N$ is indecomposable
only if $n\leq 8$. According to Proposition~2, which is applicable
(recall that $n=v_{-}+v_{+} \geq 4$), if $N$ is indecomposable, then
$S_{0}$ is neutral so that $dim \:S_{0} = 2$. Therefore, if we show that
for $n>8$ we have $dim \:S_{0} >2$, this will mean that $N$ is decomposable.

Let us complexify the source space $R^{n}$ and apply the results from \cite{1} and
\cite{2} concerning the decomposition of an $H$-normal operator in a complex
space. Lemma~1 from \cite{1} states that for an $H$-normal operator having
two distinct eigenvalues $\lambda$ and $\overline{\lambda}$ there exists a
decomposition of $C^{n}$ into a sum $C^{n}=V_{1}\dot{+}V_{2}\dot{+}V_{3}\dot{+}V_{4}$ such that
$$   N=\left( \begin{array}{cccc}
N_{1} & 0 & 0 &  0    \\
0 & N_{2} & 0 &  0    \\
0 & 0 & N_{3} &  0    \\
0 & 0 & 0 & N_{4}
\end{array} \right), \;\;
H=\left( \begin{array}{cccc}
0 & I & 0 &  0    \\
I & 0 & 0 &  0    \\
0 & 0 & H_{3} &  0    \\
0 & 0 & 0 & H_{4}
\end{array} \right),$$
where $N_{1}$, $N_{3}$ have the only eigenvalue $\lambda$, $N_{2}$, $N_{4}$
the only eigenvalue $\overline{\lambda}$, $dim V_{1}=dim V_{2}$.
It is seen that if the space $C^{n}$ is $R^{n}$ complexified, then $dim V_{3}=dim V_{4}$.

 Since ranks of the subspaces $V_{1}\dot{+}V_{2}$, $V_{3}$, $V_{4}$ are less than or
equal to 2, Theorem~1 from \cite{1} and Theorem~1 from~\cite{2} are applicable.
It follows from these theorems
that if $dim \:V_{1}, \; dim\:V_{3}>0$, then there exist at least two
linearily independent vectors $z_{1}$, $z_{2}$
such that $Nz_{1}=\lambda z_{1}$, $N^{[*]}z_{1}=\lambda z_{1}$, $Nz_{2}=\lambda z_{2}$,
$N^{[*]}z_{2}=\overline{\lambda}z_{2}$, i.e., $dim S_{0} \geq 4$.
If $dim V_{3}=0$, $n$ is equal to 4 because the subspaces $V_{1}$ and $V_{2}$
are neutral (hence $n=(2\: dim V_{1}) \leq 4 \Rightarrow n=4$). If $dim V_{1}=0$,
there appear two alternatives: $V_{3}$ and $V_{4}$ each have rank 1 or one
of these subspaces has rank 0. In the latter case either $N_{3}$ or
$N_{4}$ is decomposable for any $n$. In the former case, according to
Theorem~1~\cite{1}, $N_{3}$ ($N_{4}$) is always decomposable if
$dim V_{3}>4$ ($dim V_{4} >4$). In either case for $n>8$ there exist two linearily
independent vectors $z_{1}$, $z_{2}$ such that
$Nz_{1}=\lambda z_{1}$, $N^{[*]}z_{1}=\overline{\lambda}z_{1}$,
$Nz_{2}=\lambda z_{2}$, $N^{[*]}z_{2}=\overline{\lambda}z_{2}$. As above,
we have $dim \:S_{0} \geq 4$. Thus, if $n>8$, $N$ is decomposable, Q.E.D.

Thus, according to Proposition~2, the matrices $N$ and $H$ can be reduced to the form
\begin{equation}
N=\left( \begin{array}{ccc}
N_{1} & N_{2} & N_{3} \\
0 & N_{4} & N_{5}    \\
0 & 0 & N_{6}    \end{array}
\right), \;\;\;
H=\left( \begin{array}{ccc}
0 & 0 & I \\
0 & I & 0    \\
I & 0 & 0    \end{array}
\right),  \label{main}
\end{equation}
where
$$ N_{1}=\left( \begin{array}{cc}
\alpha & \beta \\
-\beta & \alpha \end{array}
\right), $$
$N_{6}$ is equal either to $N_{1}$ or to $N_{1}^{*}$.
The condition of the $H$-normality of $N$ is equivalent to the system
\begin{eqnarray}
N_{1}N_{6}^{*} & = & N_{6}^{*}N_{1},  \label{norm21}  \\
N_{1}N_{5}^{*}+N_{2}N_{4}^{*} & = & N_{6}^{*}N_{2}+N_{5}^{*}N_{4},  \label{norm22}  \\
N_{1}N_{3}^{*}+N_{2}N_{2}^{*}+N_{3}N_{1}^{*} & = & N_{6}^{*}N_{3}+N_{5}^{*}N_{5}+N_{3}^{*}N_{6}, \label{norm23}  \\
N_{4}N_{4}^{*} & = & N_{4}^{*}N_{4}.  \label{norm24}
\end{eqnarray}
Note that if $N_{6}=N_{1}^{*}$, then
$dim \:S_{0}''>0$ so that
it is the case $dim V_{1}>0$. It was stated before that if $dim V_{1}>0$, then
for indecomposable operators $n=4$. Therefore, for $n=4$ the submatrix $N_{6}$
can be equal to either $N_{1}$ or $N_{1}^{*}$ but for $n=6,8$ we have
$N_{6}=N_{1}$. Now let us consider the cases $n=4,6,8$ successively.

\subsubsection{$n=4$}
By the above,
$$ N= \left( \begin{array}{cc}
N_{1} & N_{3} \\
0 & N_{6} \end{array} \right)=
\left( \begin{array}{cccc}
\alpha & \beta  & a & b \\
-\beta & \alpha & c & d \\
0 & 0 & \alpha & \pm \beta \\
0 & 0 & \mp \beta & \alpha \end{array}
\right). $$

$\underline{N_{6}=N_{1}}$ Then from (\ref{norm23}) it follows that $c=-b$,
$d=a$. If $a^{2}+b^{2}=0$, i.e., $N_{3}=0$, then
$S_{0} \cap S_{1} \neq 0$, which contradicts the indecomposability of $N$.
Therefore, $a^{2}+b^{2} \neq 0$.
Taking the block diagonal transformation
$T=\sqrt[4]{a^{2}+b^{2}} I_{2} \oplus 1/\sqrt[4]{a^{2}+b^{2}} I_{2}$,
we can reduce $N$ to the form
\begin{equation}
N= \left( \begin{array}{cccc}
\alpha &  \beta &  \cos \gamma & \sin \gamma \\
-\beta & \alpha & -\sin \gamma & \cos \gamma \\
0  &     0  & \alpha    & \beta    \\
0  &     0  & -\beta    & \alpha   \end{array}
\right), \;\; 0 \leq \gamma < 2 \pi.  \label{form141}
\end{equation}

According to Proposition~3, matrix (\ref{form141}) is indecomposable.
Let us prove the $H$-unitary invariance of the parameter $\gamma$. To this
end suppose that a matrix $T$ satisfies conditions
\begin{eqnarray}
NT & = & T\tilde{N},  \label{aug_similar} \\
TT^{[*]} & = & I  \label{aug_unitar}
\end{eqnarray}
for the matrix $N$ of form
(\ref{form141}) and the matrix
$$\tilde{N}=
\left( \begin{array}{cc}
N_{1} & \widetilde{N_{3}} \\
0  &  N_{1}  \end{array} \right)=
\left( \begin{array}{cccc}
\alpha &  \beta &  \cos \tilde{\gamma} & \sin \tilde{\gamma} \\
-\beta & \alpha & -\sin \tilde{\gamma} & \cos \tilde{\gamma} \\
0  &     0  & \alpha    & \beta    \\
0  &     0  & -\beta    & \alpha   \end{array}
\right), \;\; 0 \leq \tilde{\gamma} < 2 \pi.
$$
According to Proposition~2, the matrix $T$ has the block triangular
form
$$ T= \left( \begin{array}{cc}
T_{1} & T_{2} \\
0 & T_{3} \end{array}
\right) $$
with respect to the decomposition $R^{4}=S_{0} \dot{+} S_{1}$.
The transformation $T$ is $H$-unitary iff
\begin{eqnarray}
T_{1}T_{3}^{*} & = & I,  \label{un1}  \\
T_{1}T_{2}^{*}+T_{2}T_{1}^{*} & = & 0.  \label{un2}
\end{eqnarray}
It follows from condition (\ref{aug_similar}) that $N_{1}$ and $T_{1}$
commute, therefore,
$$ T_{1}= \left( \begin{array}{cc}
t_{11} & t_{12} \\
-t_{12} & t_{11} \end{array}
\right) $$
so that from (\ref{un2}) we get
$$ T_{2}= \left( \begin{array}{cc}
t_{13} & t_{14} \\
-t_{14} & t_{13} \end{array}
\right). $$
Now, combining (\ref{un1}) and (\ref{aug_similar}), we have
$N_{1}T_{2}+N_{3}T_{1}^{*-1}=T_{1}\widetilde{N_{3}}+T_{2}N_{1}$.
But $T_{2}$ and $N_{1}$ commute (as well as $T_{1}$ and $\widetilde{N_{3}}$)
so that $N_{3}=T_{1}\widetilde{N_{3}}T_{1}^{*}=\widetilde{N_{3}}T_{1}T_{1}^{*}=(det\:T_{1})^{2}\widetilde{N_{3}}$.
Since $det N_{3}=det \widetilde{N_{3}}=1$, we have $(det T_{1})^{2}=1$
and $N_{3}=\widetilde{N_{3}}$, i.e., $\gamma=\tilde{\gamma}$, Q.E.D.

$\underline{N_{6}=N_{1}^{*}}$ Then, according to (\ref{norm23}), $c=b$.
The transformation
$$T= \left( \begin{array}{cccc}
1 & 0 & 0  & a/2\beta \\
0 & 1 & -a/2\beta & 0 \\
0 & 0 &  1  &  0 \\
0 & 0 &  0  &  1 \end{array}
\right) $$
reduces $N_{3}$ to the form
$$ N_{3}= \left( \begin{array}{cc}
0 & b' \\
b' & d' \end{array}
\right) $$
without changing the submatrices $N_{1}$ and $N_{6}$. If both $b'$ and $d'$ are
equal to zero, the condition $S_{0}\cap S_{1}=\{0\}$ fails.
Therefore, $4b'^{2}+d'^{2} \neq 0$ and we can take
the transformation
$$T= \left( \begin{array}{cccc}
\cos \phi & \sin\phi & -r \sin \phi  & r \cos \phi \\
-\sin \phi & \cos \phi & -r \cos \phi & -r \sin \phi \\
0 & 0 &  \cos \phi & \sin \phi \\
0 & 0 & -\sin \phi & \cos \phi \end{array}
\right), $$
where $\cos 2 \phi=2b'/\sqrt{d'^{2}+4b'^{2}}$,
$\sin 2 \phi=-d'/\sqrt{d'^{2}+4b'^{2}}$,
$r=d'/(4 \beta)$.
It does not change $N_{1}$ and $N_{6}$ but reduces
$N_{3}$ to the form
$$ N_{3}= \left( \begin{array}{cc}
0 & b'' \\
b'' & 0 \end{array}
\right), \; b''=\frac{1}{2}\sqrt{4b'^{2}+d'^{2}}>0. $$
If we now take $\widetilde{v_{1}}=\sqrt{b''}v_{1}$, $\widetilde{v_{2}}=\sqrt{b''}v_{2}$,
$\widetilde{v_{3}}=v_{3}/\sqrt{b''}$, $\widetilde{v_{4}}=v_{4}/\sqrt{b''}$,
then $N_{3}$ will be equal to $D_{2}$.
Thus, we have obtained the final form for the matrix $N$:
\begin{equation}
N= \left( \begin{array}{cccc}
\alpha & \beta  & 0  & 1 \\
-\beta & \alpha & 1  & 0 \\
0 & 0 & \alpha & -\beta \\
0 & 0 & \beta & \alpha \end{array}
\right). \label{form142}
\end{equation}

According to Proposition~3, matrix (\ref{form142})
is indecomposable. Forms (\ref{form141}) and (\ref{form142}) are not $H$-unitarily similar because for
matrix (\ref{form142}) the subspace $S_{0}''$ defined in Proposition~2
is nontrivial in contrast to that for (\ref{form141}).
Thus, we have proved that

{\em{If an indecomposable $H$-normal operator acts in a space $R^{4}$
of rank 2 and has 2 eigenvalues: $\alpha \pm i \beta$
($\alpha, \beta \in \Re$, $\beta > 0$),
then the pair $\{N,H\}$ is unitarily similar to one and only one of
the canonical pairs \{(\ref{lemma14.1}),(\ref{lemma14.3})\},
\{(\ref{lemma14.2}),(\ref{lemma14.3})\}.}}

\subsubsection{$n=6$}
The matrices $N$ and $H$ have form (\ref{main}) with $N_{6}=N_{1}$. Since the
submatrix $N_{4}$ is an ordinary normal matrix (condition (\ref{norm24})),
one can assume that $N_{4}=N_{1}$. Thus,
$$ N= \left( \begin{array}{ccc}
N_{1} & N_{2} & N_{3} \\
0 & N_{1} & N_{5} \\
0 & 0 & N_{1} \end{array}
\right). $$

First reduce the submatrix
$$N_{2}= \left( \begin{array}{cc}
a & b \\
c & d \end{array}
\right)
$$
to the form
\begin{equation}
N_{2}= \left( \begin{array}{cc}
0 & 0 \\
0 & 1 \end{array}
\right)  \label{n2_desired}
\end{equation}
without changing the submatrices $N_{1}=N_{4}=N_{6}$.
To this end take
\begin{equation}
T= \left( \begin{array}{ccc}
I & T_{2} & -\frac{1}{2}T_{2}T_{2}^{*} \\
0 & I & -T_{2}^{*} \\
0 & 0 & I \end{array}
\right), \label{UTF}
\end{equation}
$$ T_{2}= \left( \begin{array}{cc}
b/\beta & -a/\beta \\
0 & 0 \end{array}
\right).  $$
Then
$$N_{2}= \left( \begin{array}{cc}
0 & 0 \\
c' & d' \end{array}
\right).
$$
If both $c'$ and $d'$ are equal to zero, i.e., $N_{2}=0$, then from
condition of the $H$-normality (\ref{norm23}) it follows that $N_{5}=0$,
which contradicts the condition $S_{0} \cap S=\{0\}$.
Therefore, $c'^{2}+d'^{2} \neq 0$
and we can subject the matrix $N$ obtained to the transformation
$T=I_{2} \oplus T_{1} \oplus I_{2}$, where
$$ T_{1}= \left( \begin{array}{cc}
d'/\sqrt{c'^{2}+d'^{2}} & c'/\sqrt{c'^{2}+d'^{2}} \\
-c'/\sqrt{c'^{2}+d'^{2}} & d'/\sqrt{c'^{2}+d'^{2}} \end{array}
\right). $$
Then
$$N_{2}= \left( \begin{array}{cc}
0 & 0 \\
0 & d'' \end{array}
\right), \;\; d''=\sqrt{c'^{2}+d'^{2}} > 0.
$$
Taking $\widetilde{v_{1}}=d''v_{1}$,
$\widetilde{v_{2}}=d''v_{2}$, $\widetilde{v_{3}}=v_{3}$,
$\widetilde{v_{4}}=v_{4}$, $\widetilde{v_{5}}=v_{5}/d''$,
$\widetilde{v_{6}}=v_{6}/d''$, we obtain desired form (\ref{n2_desired})
for the submatrix $N_{2}$.

Now let us apply conditions (\ref{norm22}) and (\ref{norm23}). We get
$$N_{5}= \frac{1}{2} \left( \begin{array}{cc}
\cos \gamma + 1 & \sin \gamma \\
- \sin \gamma & \cos \gamma - 1 \end{array}
\right), \;\; 0 \leq \gamma < 2 \pi, $$
$$N_{3}= \left( \begin{array}{cc}
p & q \\
(\cos \gamma+1)/4 \beta - q & \sin \gamma/4 \beta + p \end{array}
\right).
$$

Finally, take transformation (\ref{UTF}) with
$$ T_{2}= 2p/(\cos \gamma+1)\: I_{2} \;\; \mbox{if} \;\; \gamma \neq \pi, $$
$$ T_{2}= \left( \begin{array}{cc}
0 & -q \\
q & 0 \end{array}
\right)  \;\; \mbox{if} \;\; \gamma = \pi. $$
Then $$N_{3}= \left( \begin{array}{cc}
0 & q' \\
(\cos \gamma+1)/4 \beta - q' & \sin \gamma/4 \beta \end{array}
\right) \;\; (\gamma \neq \pi), $$
$$N_{3}= p' I_{2} \;\; (\gamma = \pi). $$
As a result, we have obtained two forms:
$$ N= \left( \begin{array}{cccccc}
\alpha & \beta & 0 & 0 & 0 & r  \\
-\beta & \alpha & 0 & 1 & (\cos \gamma + 1)/4 \beta - r & \sin \gamma/4 \beta  \\
0 & 0 & \alpha & \beta & \frac{1}{2}(\cos \gamma + 1) & \frac{1}{2}\sin \gamma  \\
0 & 0 & -\beta & \alpha & -\frac{1}{2}\sin \gamma & \frac{1}{2}(\cos \gamma -1)  \\
0 & 0 & 0 & 0 & \alpha & \beta  \\
0 & 0 & 0 & 0 & -\beta & \alpha  \end{array}
\right),$$
\begin{equation}
0 \leq \gamma < 2\pi, \; \gamma \neq \pi, \label{form151}
\end{equation}
\begin{equation}
N= \left( \begin{array}{cccccc}
\alpha & \beta & 0 & 0 & r & 0  \\
-\beta & \alpha & 0 & 1 & 0 & r  \\
0 & 0 & \alpha & \beta & 0 & 0  \\
0 & 0 & -\beta & \alpha & 0 & -1  \\
0 & 0 & 0 & 0 & \alpha & \beta  \\
0 & 0 & 0 & 0 & -\beta & \alpha  \end{array}
\right).  \label{form152}
\end{equation}

According to Proposition~3, matrices (\ref{form151}) and (\ref{form152})
are indecomposable. Let us show that they are not
$H$-unitarily similar and that the numbers $r$ and $\gamma$
are $H$-unitary invariants. To this end suppose that
some $H$-unitary matrix $T$
reduces the matrix $N$ to the form $\tilde{N}$:
$$ N= \left( \begin{array}{ccc}
N_{1} & N_{2} & N_{3} \\
0 & N_{1} & N_{5} \\
0 & 0 & N_{1} \end{array}
\right), \;\;
\tilde{N}= \left( \begin{array}{ccc}
N_{1} & N_{2} & \widetilde{N_{3}} \\
0 & N_{1} & \widetilde{N_{5}} \\
0 & 0 & N_{1} \end{array}
\right), $$
where
$$ N_{1}= \left( \begin{array}{cc}
\alpha & \beta \\
-\beta & \alpha \end{array}
\right), \;\;
N_{2}= \left( \begin{array}{cc}
0 & 0 \\
0 & 1 \end{array}
\right), $$
$$N_{5}= \frac{1}{2} \left( \begin{array}{cc}
\cos \gamma + 1 & \sin \gamma \\
- \sin \gamma & \cos \gamma - 1 \end{array}
\right), \;\; 0 \leq \gamma < 2 \pi,  $$
$$\widetilde{N_{5}}= \frac{1}{2} \left( \begin{array}{cc}
\cos \tilde{\gamma} + 1 & \sin \tilde{\gamma} \\
- \sin \tilde{\gamma} & \cos \tilde{\gamma} - 1 \end{array}
\right), \;\; 0 \leq \tilde{\gamma} < 2 \pi.  $$
Then, according to Proposition~2, $T$ has the block triangular form
$$T= \left( \begin{array}{ccc}
T_{1} & T_{2} & T_{3} \\
0 & T_{4} & T_{5} \\
0 & 0 & T_{6} \end{array}
\right)  $$
with respect to the decomposition $R^{6}=S_{0}\dot{+}S\dot{+}S_{1}$.
It follows from condition (\ref{aug_similar}) $NT=T\tilde{N}$ that
$$ T_{1}=T_{4}=T_{6}= \left( \begin{array}{cc}
\cos \phi & \sin \phi \\
-\sin \phi & \cos \phi \end{array}
\right), \;\;
T_{2}= \left( \begin{array}{cc}
t_{13} & t_{14} \\
-t_{14} & t_{13}+\frac{\sin \phi}{\beta} \end{array}
\right). $$
Condition (\ref{aug_unitar}) $TT^{[*]}=I$ implies
$T_{1}T_{5}^{*}+T_{2}T_{4}^{*}=0$, hence
$$ T_{5}=-T_{1}T_{2}^{*}T_{1}=
\left( \begin{array}{cc}
t_{35} & t_{36} \\
-t_{36} & t_{35}-\frac{\sin \phi}{\beta} \end{array}
\right), $$
where
\begin{eqnarray*}
t_{35} & = & -t_{13}\cos 2\phi-t_{14} \sin 2\phi + \frac{\sin^{3} \phi}{\beta}, \\
t_{36} & = & -t_{13}\sin 2 \phi +t_{14} \cos 2\phi - \frac{\cos \phi \sin^{2} \phi}{\beta}.
\end{eqnarray*}
Substituting the expressions for $T_{4}$, $T_{5}$, $T_{6}$ in the
formula $N_{1}T_{5}+N_{5}T_{6}=T_{4}\widetilde{N_{5}}+T_{5}N_{1}$, which
follows from (\ref{aug_similar}), we obtain: $N_{5}=\widetilde{N_{5}}$.
Therefore, forms (\ref{form151}) and (\ref{form152}) are not $H$-unitarily
similar and the parameter $\gamma$ is an $H$-unitary invariant.

Now let us check the $H$-unitary invariance of $r$ for matrix (\ref{form151}).
To this end suppose that
$$N_{3}= \left( \begin{array}{cc}
0 & r \\
(\cos \gamma+1)/4 \beta - r & \sin \gamma/4 \beta \end{array}
\right), $$
$$\widetilde{N_{3}}= \left( \begin{array}{cc}
0 & \tilde{r} \\
(\cos \gamma+1)/4 \beta - \tilde{r} & \sin \gamma/4 \beta \end{array}
\right),$$
$0 \leq \gamma < 2 \pi$, $\gamma \neq \pi$.
It follows from (\ref{aug_unitar}) that $T_{1}T_{3}^{*}=-\frac{1}{2}T_{2}T_{2}^{*}+X$,
where $X$ is an antisymmetric matrix, therefore,
$$T_{3}= \left( \begin{array}{cc}
t_{15} & t_{16} \\
t_{25} & t_{26} \end{array}
\right) -
\left( \begin{array}{cc}
x \sin \phi & -x \cos \phi \\
x \cos \phi & x \sin \phi \end{array}
\right), $$
where
\begin{eqnarray*}
2t_{15} & = & -(t_{13}^{2}+t_{14}^{2})\cos \phi +t_{14} \sin^{2} \phi /\beta, \\
2t_{16} & = & -(t_{13}^{2}+t_{14}^{2})\sin \phi -t_{14}\sin \phi \cos \phi /\beta, \\
2t_{25} & = & -t_{14}\sin \phi \cos \phi/\beta + ((t_{13}+\sin \phi/\beta)^{2}+t_{14}^{2})\sin \phi, \\
2t_{26} & = & -t_{14} \sin^{2} \phi /\beta - ((t_{13}+ \sin \phi/\beta)^{2}+t_{14}^{2})\cos \phi.
\end{eqnarray*}
Since $N_{1}T_{3}+N_{2}T_{5}+N_{3}T_{6}=T_{1}\widetilde{N_{3}}+T_{2}N_{5}+
T_{3}N_{1}$ (condition (\ref{aug_similar})),
$\widetilde{N_{3}}=T_{1}^{*}(N_{1}T_{3}-T_{3}N_{1}+N_{2}T_{5}-T_{2}N_{5}+N_{3}T_{6})$.
Substituting the expressions
for $T_{2}$, $T_{3}$, $T_{5}$, $T_{6}$ in this formula, we obtain:
\begin{eqnarray}
a_{1}t_{13}+a_{2}t_{14}+a_{3} & = & 0,  \label{sss1} \\
b_{1}t_{13}+b_{2}t_{14}+b_{3} & = & \tilde{r}-r, \label{sss2}
\end{eqnarray}
where
\begin{eqnarray*}
a_{1} & = & -\frac{1}{2}(\cos(\phi - \gamma)+\cos\phi), \\
a_{2} & = & -\frac{1}{2}(\sin(\phi-\gamma)+\sin\phi), \\
a_{3} & = & -\frac{1}{4 \beta} \sin \phi (\cos (\phi - \gamma)+\cos \phi), \\
b_{1} & = & \frac{1}{2}(\sin(\phi -\gamma)-\sin \phi), \\
b_{2} & = & -\frac{1}{2}(\cos(\phi -\gamma)-\cos \phi), \\
b_{3} & = & \frac{1}{4 \beta} \sin \phi (\sin (\phi-\gamma)-\sin \phi).
\end{eqnarray*}
Since the left hand sides of equations (\ref{sss1}) - (\ref{sss2}) are proportional
and the coefficients of $t_{13}$ and of $t_{14}$ in (\ref{sss1}) are not equal
to zero simultaneously, condition (\ref{sss1}) implies $\tilde{r}=r$.
Therefore, $r$ is an $H$-unitary invariant. The proof of the invariance of $r$
for matrix (\ref{form152}) is analogous.

Thus, we have proved that

{\em{If an indecomposable $H$-normal operator acts in a space $R^{6}$
of rank 2 and has 2 eigenvalues: $\alpha \pm i \beta$
($\alpha, \beta \in \Re$, $\beta > 0$),
then the pair $\{N,H\}$ is unitarily similar to one and only one of
the canonical pairs \{(\ref{lemma15.1}),(\ref{lemma15.3})\},
\{(\ref{lemma15.2}),(\ref{lemma15.3})\}.}}

\subsubsection{$n=8$}
The matrices $N$ and $H$ have form (\ref{main}), $N_{6}$ being equal to $N_{1}$:
$$ N= \left( \begin{array}{ccc}
N_{1} & N_{2} & N_{3}  \\
0 & N_{4} & N_{5} \\
0 & 0 & N_{1} \end{array}
\right). $$
Since $N_{4}$ is an ordinary normal matrix (condition (\ref{norm24})), it can be assumed
that $N_{4}=N_{1} \oplus N_{1}$.

Having these equalities in mind, we reduce the submatrix
$$N_{2}= \left( \begin{array}{cccc}
a & b & c & d \\
e & f & g & h \end{array}
\right)
$$
to the form
\begin{equation}
N_{2}= \left( \begin{array}{cccc}
0 & 0 & 0 & 0 \\
0 & 1 & 0 & 0 \end{array}
\right)  \label{n2_new}
\end{equation}
without changing the submatrices $N_{1}$, $N_{4}$, and $N_{6}=N_{1}$.
To this end take transformation (\ref{UTF}) with
$$ T_{2}= \left( \begin{array}{cccc}
b/\beta & -a/\beta & d/\beta & -c/\beta \\
0       &  0       & 0       & 0         \end{array}
\right).  $$
Then
$$N_{2}= \left( \begin{array}{cccc}
0 & 0 & 0 & 0 \\
e' & f' & g' & h' \end{array}
\right).
$$
Now subject the obtained matrix $N$ to the
transformation $T=I_{2}\oplus T_{1}' \oplus T_{1}'' \oplus I_{2}$,
where
$$ T_{1}'=\left( \begin{array}{cc}
f'/\sqrt{e'^{2}+f'^{2}} & e'/\sqrt{e'^{2}+f'^{2}} \\
-e'/\sqrt{e'^{2}+f'^{2}} & f'/\sqrt{e'^{2}+f'^{2}} \end{array}
\right) \;\; \mbox{if} \; e'^{2}+f'^{2}>0, $$
$$ T_{1}'=I_{2} \;\; \mbox{if} \; e'=f'=0, $$
$$ T_{1}''= \left( \begin{array}{cc}
h'/\sqrt{g'^{2}+h'^{2}} & g'/\sqrt{g'^{2}+h'^{2}} \\
-g'/\sqrt{g'^{2}+h'^{2}} & h'/\sqrt{g'^{2}+h'^{2}} \end{array}
\right) \;\; \mbox{if} \; g'^{2}+h'^{2}>0,  $$
$$ T_{1}''=I_{2} \;\; \mbox{if} \; g'=h'=0. $$
We get
$$N_{2}= \left( \begin{array}{cccc}
0 & 0 & 0 & 0 \\
0 & f''& 0 & h'' \end{array}
\right), \; f''=\sqrt{e'^{2}+f'^{2}} \geq 0, \; h''=\sqrt{g'^{2}+h'^{2}} \geq 0.
$$
If $f''+h''=0$, i.e., $N_{2}=0$, from condition (\ref{norm23}) it follows that
$N_{5}=0$, which is impossible because $S_{0}\cap S=\{0\}$.
Therefore, $f''+h''>0$. Without loss of generality it can be assumed that
$f'' \neq 0$ (otherwise one can take $\widetilde{v_{3}}=v_{5}$,
$\widetilde{v_{4}}=v_{6}$, $\widetilde{v_{5}}=v_{3}$,
$\widetilde{v_{6}}=v_{4}$). Therefore, we can assume $f''=1$, taking
$\widetilde{v_{1}}=f''v_{1}$, $\widetilde{v_{2}}=f''v_{2}$,
$\widetilde{v_{7}}=v_{7}/f''$, $\widetilde{v_{8}}=v_{8}/f''$.
Keeping in mind that $f''=1$, take the transformation
$$ T= T_{1} \oplus \left( \begin{array}{cc}
1/\sqrt{1+h''^{2}} I_{2} & -h''/\sqrt{1+h''^{2}} I_{2} \\
h''/\sqrt{1+h''^{2}} I_{2} & 1/\sqrt{1+h''^{2}} I_{2}  \end{array}
\right) \oplus T_{1}^{*-1}, $$
where $T_{1}=\sqrt{1+h''^{2}} I_{2}$.
Then we obtain desired form (\ref{n2_new}) for the submatrix $N_{2}$.

Condition (\ref{norm22}) implies
$$N_{5}= \left( \begin{array}{cc}
{*} & {*} \\
{*} & {*} \\
p  &  q  \\
- q &  p  \end{array}
\right). $$
Since the case $p=q=0$ is impossible (the condition $S_{0}\cap S=\{0\}$),
we have $p^{2}+q^{2} > 0$. The transformation
$T=I_{2} \oplus I_{2} \oplus T_{1} \oplus I_{2}$, where
$$ T_{1}= \left( \begin{array}{cc}
p/\sqrt{p^{2}+q^{2}}  & q/\sqrt{p^{2}+q^{2}} \\
-q/\sqrt{p^{2}+q^{2}} & p/\sqrt{p^{2}+q^{2}} \end{array}
\right), $$
reduces $N_{5}$ to the form
$$N_{5}= \left( \begin{array}{cc}
{*} & {*} \\
{*} & {*} \\
p' &  0  \\
0  &  p' \end{array}
\right), \; p'=\sqrt{p^{2}+q^{2}} >0, $$
retaining the submatrices $N_{1}$, $N_{2}$, $N_{4}$, and $N_{6}$.
It follows from conditions of the $H$-normality (\ref{norm22}) and
(\ref{norm23}) that
$$N_{5}= \left( \begin{array}{cc}
\sin^{2} \gamma & \sin \gamma \cos \gamma \cos \delta \\
-\sin \gamma \cos \gamma \cos \delta & -\cos^{2} \gamma  \\
\sin \gamma \cos \gamma \sin \delta &  0  \\
0  &  \sin \gamma \cos \gamma \sin \delta \end{array}
\right), \; 0< \gamma < \pi/2, \; 0< \delta < \pi,  $$
$$N_{3}= \left( \begin{array}{cc}
s & t \\
\sin^{2} \gamma/2\beta - t & \sin \gamma \cos \gamma \cos \delta / 2 \beta + s  \end{array}
\right). $$

At last, take transformation (\ref{UTF}), where
$$ T_{2}= \left( \begin{array}{cccc}
0 & 0 & s/(\sin \gamma \cos \gamma \sin \delta) & t/(\sin \gamma \cos \gamma \sin \delta) \\
0 & 0 & -t/(\sin \gamma \cos \gamma \sin \delta) & s/(\sin \gamma \cos \gamma \sin \delta) \end{array} \right), $$
and reduce $N$ to the final form:
$$ N= \left( \begin{array}{cccccccc}
\alpha & \beta & 0 & 0  & 0 & 0 & 0 & 0  \\
-\beta & \alpha & 0 & 1  & 0 & 0 & \sin^{2}\gamma/2\beta & \sin \gamma \cos \gamma \cos \delta /2 \beta  \\
0 & 0 & \alpha & \beta  & 0 & 0 & \sin^{2} \gamma & \sin \gamma \cos \gamma \cos \delta  \\
0 & 0 & -\beta & \alpha & 0 & 0 & -\sin \gamma \cos \gamma \cos \delta & -\cos^{2} \gamma  \\
0 & 0 & 0 & 0 & \alpha & \beta & \sin \gamma \cos \gamma \sin \delta & 0   \\
0 & 0 & 0 & 0 & -\beta & \alpha & 0 & \sin \gamma \cos \gamma \sin \delta  \\
0 & 0 & 0 & 0 & 0 & 0 & \alpha & \beta   \\
0 & 0 & 0 & 0 & 0 & 0 & -\beta & \alpha  \end{array}
\right) $$
\begin{equation}
0<\gamma<\pi/2, \; 0<\delta <\pi. \label{form161}
\end{equation}

Due to Proposition~3 the matrix obtained is indecomposable.
Let us check the $H$-unitary invariance of the parameters
$\gamma$ and $\delta$. Suppose some $H$-unitary matrix $T$
reduces the matrix $N$ to the form $\tilde{N}$:
$$ N= \left( \begin{array}{ccc}
N_{1} & N_{2} & N_{3} \\
0 & N_{4} & N_{5} \\
0 & 0 & N_{1} \end{array}
\right), \;\;
\tilde{N}= \left( \begin{array}{ccc}
N_{1} & N_{2} & \widetilde{N_{3}} \\
0 & N_{4} & \widetilde{N_{5}} \\
0 & 0 & N_{1} \end{array}
\right), $$
where
$$ N_{1}= \left( \begin{array}{cc}
\alpha & \beta \\
-\beta & \alpha \end{array}
\right), \;\;
N_{4}=N_{1} \oplus N_{1}, \;\;
N_{2}= \left( \begin{array}{cccc}
0 & 0 & 0 & 0 \\
0 & 1 & 0 & 0 \end{array}
\right), $$
$$N_{5}= \left( \begin{array}{cc}
\sin^{2} \gamma & \sin \gamma \cos \gamma \cos \delta \\
-\sin \gamma \cos \gamma \cos \delta & -\cos^{2} \gamma  \\
\sin \gamma \cos \gamma \sin \delta &  0  \\
0  &  \sin \gamma \cos \gamma \sin \delta \end{array}
\right), $$
$$ \widetilde{N_{5}}= \left( \begin{array}{cc}
\sin^{2} \tilde{\gamma} & \sin \tilde{\gamma} \cos \tilde{\gamma} \cos \tilde{\delta} \\
-\sin \tilde{\gamma} \cos \tilde{\gamma} \cos \tilde{\delta} & -\cos^{2} \tilde{\gamma}  \\
\sin \tilde{\gamma} \cos \tilde{\gamma} \sin \tilde{\delta} &  0  \\
0  &  \sin \tilde{\gamma} \cos \tilde{\gamma} \sin \tilde{\delta} \end{array}
\right),  $$
$$ 0< \gamma, \tilde{\gamma} < \pi/2, \; 0< \delta, \tilde{\delta} < \pi. $$

Then, according to Proposition~2, $T$ has the block triangular form
$$T= \left( \begin{array}{ccc}
T_{1} & T_{2} & T_{3} \\
0 & T_{4} & T_{5}  \\
0 & 0  &  T_{6} \end{array} \right) $$
with respect to the decomposition $R^{8}=S_{0}\dot{+}S\dot{+}S_{1}$.
Combining condition (\ref{aug_similar}) $NT=T\tilde{N}$ and
(\ref{aug_unitar}) $TT^{[*]}=I$, we get
$$ T_{1}=T_{6}= \left( \begin{array}{cc}
\cos \phi & \sin \phi \\
-\sin \phi & \cos \phi \end{array}
\right), \; \;
T_{4}=T_{1}\oplus
\left( \begin{array}{cc}
\cos \psi & \sin \psi \\
-\sin \psi & \cos \psi \end{array}
\right), $$
$$ T_{2}= \left( \begin{array}{cccc}
t_{13} & t_{14} & t_{15} & t_{16} \\
-t_{14} & t_{13}+\frac{\sin \phi}{\beta} & -t_{16} & t_{15} \end{array}
\right), \; \;
T_{5}= \left( \begin{array}{cc}
t_{37} & t_{38} \\
-t_{38} & t_{37}-\frac{\sin \phi}{\beta} \\
t_{57} & t_{58} \\
-t_{58} & t_{57} \end{array} \right), $$
where
\begin{eqnarray*}
t_{37} & = & -t_{13}\cos 2 \phi -t_{14}\sin 2 \phi +\frac{\sin^{3}\phi}{\beta}, \\
t_{38} & = & -t_{13}\sin 2 \phi +t_{14}\cos 2 \phi - \frac{\cos \phi \sin^{2} \phi}{\beta}, \\
t_{57} & = & -t_{15}\cos (\phi+\psi) - t_{16}\sin (\phi+\psi), \\
t_{58} & = & -t_{15}\sin (\phi+\psi) + t_{16}\cos (\phi+\psi).
\end{eqnarray*}
Substituting the expressions for  $T_{4}$, $T_{5}$, $T_{6}$
in the formula $N_{4}T_{5}+N_{5}T_{6}=T_{4}\widetilde{N_{5}}+T_{5}N_{1}$
which follows from  (\ref{aug_similar}), we
obtain
$$ \widetilde{N_{5}}= \left( \begin{array}{cc}
\sin^{2} \gamma & \sin \gamma \cos \gamma \cos \delta \\
-\sin \gamma \cos \gamma \cos \delta & -\cos^{2} \gamma \\
\sin \gamma \cos \gamma \sin \delta \cos (\phi-\psi) & \sin \gamma \\cos \gamma \sin \delta \sin (\phi-\psi) \\
-\sin \gamma \cos \gamma \sin \delta \sin (\phi-\psi) & \sin \gamma \cos \gamma \sin \delta \cos (\phi-\psi) \end{array} \right), $$
hence $\phi=\psi$, hence $\gamma=\tilde{\gamma}$, $\delta=\tilde{\delta}$.

Thus, we have proved that

{\em{If an indecomposable $H$-normal operator acts in a space $R^{8}$
of rank 2 and has 2 eigenvalues: $\alpha \pm i \beta$
($\alpha, \beta \in \Re$, $\beta > 0$),
then the pair $\{N,H\}$ is unitarily similar to the canonical
pair \{(\ref{lemma16.1}),(\ref{lemma16.2})\}.}}

We have considered all the possible alternatives for an indecomposable
operator $N$ and have obtained the canonical forms for each case. Thus,
we have proved Theorem~2.


\begin{thebibliography}{1}
\bibitem{1} I. Gohberg, B. Reichstein, {\em{On classification of
Normal Matrices in an Indefinite Scalar Product, \/}} Integral
Equations and Operator Theory, 13 (1990), 364-394.

\bibitem{2} O.V.Holtz, V.A.Strauss, {\em{Classification of Normal 
Operators in Spaces with Indefinite Scalar Product of Rank 2, \/}} 
Linear Algebra Appl., 241--3 (1996), 455-517.


\end{thebibliography}
\end{document}